%% file: main.tex
\documentclass[review]{elsarticle}

\usepackage{lineno,hyperref}
\usepackage{amsmath}
\usepackage{breqn}
\usepackage{calc}  
\usepackage{enumitem}  
\usepackage{amsfonts}
\usepackage{mathrsfs}
\usepackage{algpseudocode}
\usepackage{multirow}
\usepackage{bigstrut}
\usepackage{booktabs}
\usepackage{rotating}
\usepackage{hyperref}
\usepackage{bm}
\usepackage[top=3.5cm,bottom=3.5cm,left=3.5cm,right=3.5cm]{geometry}
\usepackage{caption}
\usepackage{subcaption}
\usepackage{mathtools}
\usepackage{algorithm2e}

\input{mac}

 \usepackage{caption}

\biboptions{authoryear}

\makeatletter
\def\ps@pprintTitle{%
 \let\@oddhead\@empty
 \let\@evenhead\@empty
 \def\@oddfoot{}%
 \let\@evenfoot\@oddfoot}
\makeatother
\begin{document}

\begin{frontmatter}

\title{Spatial Pricing of Ride-sourcing Services in a Congested Transportation Network}



\author[mymainaddress]{Fatima Afifah}
\author[mymainaddress]{Zhaomiao Guo\corref{mycorrespondingauthor}}
\cortext[mycorrespondingauthor]{(Corresponding Author) Assistant Professor, 4353 Scorpius Street, Orlando, Florida 32816-0120, Phone: 407-823-6215, Email: guo@ucf.edu}

\address[mymainaddress]{Department of Civil, Environmental and Construction Engineering\\
Resilient, Intelligent and Sustainable Energy Systems Cluster\\ 
	University of Central Florida, FL 32766
}

\begin{abstract}
	 We investigate the impacts of spatial pricing for ride-sourcing services in a Stackelberg framework considering traffic congestion. In the lower level, we use combined distribution and assignment approaches to explicitly capture the interactions between drivers' relocation, riders' mode choice, and all travelers' routing decisions. In the upper level, a monopolistic transportation network company (TNC) determines spatial pricing strategies to minimize imbalance in a two-sided markets. We show the existence of the optimal pricing strategies for locational imbalance minimization, and propose effective algorithms with reliable convergence properties. Furthermore, the optimal pricing is unique, and can be solved in a convex reformulation when matching time can be ignored. We conduct numerical experiments on different scales of transportation networks with different TNC objectives to generate policy insights on how spatial pricing could impact transportation systems. 
\end{abstract}

\end{frontmatter}


\section{Introduction}

Since the first introduction of UBER services in 2010, ride-sourcing has become an encouraging new form of mobility. Transportation network companies (TNCs), also known as ride-sourcing companies, serve as both platform providers and operators, promising to increase reliability and efficiency of point-to-point transportation. Dynamic pricing\footnote{E.g., Uber surge pricing, Lyft prime time charges, and Didi Chuxing's dyanamic pricing.} is one of TNCs' key operational techniques to improve locational supply-demand balance in real-time. Conceptually, dynamic pricing sets a surge multiplier (SM) for each predefined geographic area, which dynamically raises (or reduces) the base trip fare when locational demand is higher (or lower) than available supply until demand and supply are balanced \footnote{In reality, prices for rides are dynamically calculated based on a variety of factors including route, time of day, ride type, number of available drivers, current demand for rides, and any local fees or surcharges. While actual pricing algorithm is unknown to public, \cite{chen2015peeking} reverse engineer the pricing strategies using high-resolution UBER data.}. While conceptually appealing, dynamic pricing receives critiques over the past years due to the opaque of pricing algorithm and concerns that TNCs may surge unnecessarily high or more frequently to exploit customers in an oligopolistic setting \citep{zha2018surge}. 

Determining optimal pricing strategies may be challenging due to spatio-temporal couplings and potential conflict of interests between the society and TNCs. On one hand, dynamic pricing tries to balance supply and demand both in time and in space, but the most research attention has been paid to the temporal aspect \citep{bimpikis2019spatial}. Because pricing information (i.e., SM, base trip fare, etc.) are revealed to both demand and supply markets, which are interconnected through a transportation network, prices at one location affect services at all other locations potentially. In practice, TNCs also benefit from using spatial prices to account for long-term predictable unbalanced demand patterns\citep{bimpikis2019spatial}. Studying dynamic pricing from a network perspective is beneficial for both societal welfare and TNCs' mission. On the other hand, while TNCs are typically profit-driven, optimal prices for TNCs themselves may not be optimal for the whole society, which includes riders, drivers and other transportation users. Better understanding the impact of different TNCs' objectives on system welfare is critical to generate planning and regulation insights for policy makers.

In this paper, we focus on addressing the spatial effects of dynamic pricing considering traffic congestion. We investigate the impacts of spatial pricing \footnote{Spatial pricing is the average pricing of ride-sourcing services originated from predefined zones, such as census tracts or 1mi by 1mi blocks. We note that spatial pricing and surge multipliers are closely related. If we ignored the temporal aspect of surge multipliers, the ratio between spatial pricing and average trip fares in normal states should be close to surge multipliers.} for ride-sourcing services in an interconnected transportation network using a Stackelberg framework. In the lower level, we explicitly capture the interactions between drivers' relocation, riders' mode choice, and all travelers' routing decisions. In the upper level, a monopolistic TNC determines spatial pricing strategies to fulfill its own objective in a two-sided markets. The main contribution of this work is threefold. First, we explicitly consider the relocation of drivers and transportation congestion, which not only more accurately captures the ride-sourcing behaviors, but also allows for systematic analyses of the impacts of spatial pricing on transportation efficiency. Second, while optimal spatial pricing problems are challenging to solve due to a bi-level structure, we propose effective computational algorithms, with rigorous analyses on the existence and convergence of the optimal solutions. Third, when matching time can be ignored, we reformulate the optimal spatial pricing problem as a single-level convex optimization, which significantly improves the computational efficiency and guarantees a unique global optimal.

The remainder of this paper is organized as follows. Section 2 reviews the relevant literature on dynamic pricing in ride-sourcing context. Section 3 and section 4 propose the formulations and algorithms, respectively. We test the models and solution approaches using difference scales of transportation systems in section 5. Section 6 concludes the paper with main research findings and future research directions.

\section{Literature Review}

According to \cite{sae2018taxonomy}, ride-sourcing is defined as ``prearranged and on-demand transportation services for compensation in which drivers and passengers connect via digital applications''. This definition includes a broad range of mobility modes. In this paper, we focus on TNCs for automobiles.

As one of the key operational strategies, ride-sourcing pricing has attracted increasing attention from various research community, including transportation (e.g., \citep{zha2018surge, zha2018geometric}), economics (e.g., \citep{chen2016dynamic, bimpikis2019spatial}), management science (e.g., \citep{guda2019your}) and computer science (e.g., \citep{chen2015peeking}). In addition, some studies on ride-sourcing modeling are closely related to research on taxi (e.g., \citep{yang1998network, lagos2000alternative,  yang2010equilibria, yang2011equilibrium, he2018pricing}) and two-sided marketplace (e.g., \cite{rochet2006two, rysman2009economics}). Given an increasingly large body of literature, we focus on the pricing issues in the context of ride-sourcing from transportation perspective, with references to other domains as needed.


The impacts of dynamics pricing on ride-sourcing driver supply are still inconclusive. The first conjecture, denoted as income-target theory, is that drivers quit driving once a daily income target is reached. This theory is pioneered by \cite{camerer1997labor}, who find negative wage elasticity of cabdriver supply in New York City. In contrast, the second conjecture, denoted as the neoclassical theory\citep{farber2005tomorrow}, suggests that cabdrivers are greedy for higher income, and cabdrivers supply has a positive wage rate elasticity. \cite{farber2005tomorrow} also points out that the differences between income-target theory and neoclassical theory are in the conception and measurement of the daily wage rate.  \cite{chen2016dynamic} provide empirical evidence for neoclassical theory using high-resolution UBER data. The coexistence of these two competing theories may explain the empirical evidence that the effect of dynamic pricing on UBER drivers supply is limited\citep{chen2015peeking}.  \cite{zha2018surge} propose a bi-level modeling framework, with lower level equilibrium considering both income-target theory and neoclassical theory, and find that compared to static pricing, surge pricing leads to higher revenue for both the platform and drivers, but lower customer surplus during highly surged periods. However, the majority of these studies focus on temporal dimension of dynamic pricing without considering spatial linkages or driver relocations. In this paper, due to the inconclusive impacts of dynamic pricing on driver supply, we assume that total driver supply is inelastic to spatial pricing, but drivers will respond to spatial pricing by choosing their service location.

Another school of literature investigate the potential value of dynamic pricing when market condition is stochastic\citep{banerjee2015pricing, cachon2017role, castillo2017surge, taylor2018demand, gurvich2019operations}. \cite{banerjee2015pricing} adopt a queueing-theoretical approach and find that dynamic pricing only outperform static pricing when TNCs have imperfect knowledge of system parameters. Although network model is presented, the pricing policies only depend on number of idle drivers at individual locations. \cite{castillo2017surge} consider the problem of ``wild goose chase'' when price is set too low, and find that dynamic pricing can avoid this issue and maintaining prices within a reasonable range. In contrast, in a deterministic setting, \cite{zha2018geometric} adopt a geometric matching framework and find that dynamic pricing may still set price higher than the societal efficient level. However, idle vehicle relocation behaviors is not considered in \citep{zha2018geometric}.

\cite{bimpikis2019spatial} are among the first to consider the spatial dimension of ride-sourcing pricing, by modeling whether and where drivers will join the platform, and where to relocate themselves when they are idle. \cite{bimpikis2019spatial} find that spatial pricing contributes to balance supply and demand, and more balanced demand patterns lead to higher profits and higher consumer surplus, simultaneously. However, the equilibrium formulation in \citep{bimpikis2019spatial} is indeterminate and the conclusion is, therefore, best seen as an ideal bound\citep{zha2018geometric}. In addition, since \cite{bimpikis2019spatial} mainly focus on the equilibrium of ride-sourcing services rather than traffic equilibrium, transportation topology and congestion effects are ignored. \cite{nie2017can} has indicated the concerns of ride-sourcing services on worsening the traffic congestion. A more recent working paper by \cite{ban2018general} using general equilibrium model to formulate the multi-agent interactions in an e-hailing system, has explicitly modeled traffic congestion under path independent free-flow travel time assumption. However, pricing are treated as exogenous parameters and not a focus in \citep{ban2018general}. We aim to fill this gap by proposing a mathematical tractable modeling framework for spatial pricing considering traffic equilibrium in a congested transportation network.

\section{Methodologies}\label{sec:meth}

\subsection{Problem Settings and Key Assumptions}
We assume that only one TNC serves in our study area, who aims to optimally decide the spatial prices to minimize demand and supply imbalance. We will consider different objectives in Section \ref{sec:nume} to illustrate the impacts of different TNC objectives on spatial pricing strategies. Due to the uncoordinated behaviors of a large number of drivers and riders, individual driver and rider do not have market power, and only respond to the spatial prices set by the TNC. We assume drivers make their relocation/pickup decisions based on locational attractiveness, relocation time, waiting time, and locational prices. Locational attractiveness could include factors such as safety, rider density, easiness of getting returning trips. The modeling framework is flexible to incorporate more explanatory factors if needed.  While drivers have flexibility to relocate themselves to other locations, travelers are typically not able to change their travel origin. We assume that the travel demand not requesting TNC services will shift to car driving. The interactions of drivers, riders, and other passenger vehicles will endogeneously determine the mobility of the transportation system. Notice that although trips originated from one location may have different trip distances with different prices, we assume that drivers and travelers will take the average prices originating from one location into their decision making processes. This assumption is consistent with the current practice that SMs are based on trip origins only. Heterogeneous demand, for example, trips with longer travel distance are more sensitive to spatial prices, can be incorporated into our proposed modeling framework without fundamental changes on  modeling and computational strategies. 

In the remaining of this section, we present the mathematical models for two cases: ignoring and considering matching time.  Ignoring matching time, which is significantly easier to solve (see in Section \ref{sec:solu_meth}),  is appropriate especially when congestion is severe, and drivers and riders are balanced in one region (therefore waiting time is small compared to relocation/travel time). In addition, some riders may not be sensitive to waiting time as they can request rides earlier than their scheduled departure time. We note that ignored matching time is a special case of considering matching time case.

\subsection{Mathematical Modeling: Ignoring Matching Time}

For long-term transportation planning problem with congestion, the relocation time can significantly exceed the matching time. In this section, we ignore the impacts of matching time on the behaviors of drivers and riders. This assumption will be relaxed in Section \ref{sec:modeling_matching}.

\subsubsection{Traffic Modeling}
Denote a transportation network by a directed graph $\mathcal{G} = (\mathcal{N}, \mathcal{A})$ , where $\mathcal{N}$  is the set of nodes (indexed by $n$) and  $\mathcal{A}$ is the set of links (indexed by $a$).  The service area we consider is a metropolitan area where ride-sourcing demand and supply are aggregated in predefined zones (e.g., census tracts or 1 mile by 1 mile blocks), represented by nodes.  A link may represent a road section that connects two nodes. Given driver node set $\mathcal{R}$ ($\subset \mathcal{N}$) and rider node set $\mathcal{S}$ ($\subset \mathcal{N}$), drivers will relocate from $r$ ($\in \mathcal{R}$) to $s$ ($\in \mathcal{S}$) to pick up riders. Note that the relocation of drivers could be due to accepted matchings or voluntary relocations. Also note that if a driver decides to wait at his/her current location, or has accepted trip in the same area, $r = s$.  We assume that driver relocation decisions following multinomial logistic models, with explanatory variables of locational attractiveness $\beta_{0,s}$, relocation time $t_{rs}$, and locational prices $\rho_s$. The deterministic components of drivers' utility function are given in (\ref{eq:utility}).

\begin{equation}
\label{eq:utility}
U_{rs} = \beta_{0,s} -\beta_1 t_{rs} + \beta_2 \rho_{s}
\end{equation}
where: 
\begin{description}[leftmargin=!,labelwidth=\widthof{12345}]
	\itemsep-0.7em 
	\item[$U_{rs}$]: deterministic component of utility measure for a driver going from $r$ to $s$;
	\item[$\beta$] : utility function parameters (model input);
	\item[$t_{rs}$] : equilibrium travel time from $r$ to $s$; 
	\item[$\rho_{s}$] : locational price at $s$.
\end{description}

While drivers have flexibility to relocate themselves to other locations, we assume riders are not able to change their travel origins. At each location $s$ ($\in \mathcal{S}$), we  assume a linear ride-sourcing demand function, as defined in (\ref{eq:demand}). One can adopt more sophisticated demand models, such as binary logit model in Section \ref{sec:modeling_matching}, to describe the ride-sourcing demand as a (decreasing) function of prices, which will lead to a similar single-level convex reformulation as illustrated in Section \ref{sec:solu_approach_ignore}.

\begin{equation}
\label{eq:demand}
d_{s} = D_s - b_s\rho_{s}
\end{equation}
where: 
\begin{description}[leftmargin=!,labelwidth=\widthof{12345}]
	\itemsep-0.7em 
	\item[$d_{s}$]: number of riders requested service at location $s$;
	\item[$D_s, b_s$] : coefficients of demand function (model input).
\end{description}

We assume that the total vehicle travel demand, including ride-sourcing travel and passenger vehicle travel, at location $s$ is given and denoted as $\bar{D}_s$. After been picked up by drivers, ride-sourcing riders, as well as other passenger vehicles, will travel at routes that minimize their own travel time. Since $t_{rs}$ is endogenously determined by traffic equilibrium, we extend combined distribution and assignment (CDA) model \citep{wilson1969use} to describe the interactions between drivers' relocation and the routing of all ride-sourcing vehicles and conventional passenger vehicles. The model is described in  (\ref{model:lower}).

\begin{subequations}
\begin{align}
& & & \underset{\bm{\hat{v}},\bm{\check{v}}, \bm{q} \in \reals_+}{\text{minimize}}
& & \sum_{a \in \mathcal{A}} \int_{0}^{v_a} t_a(v_a) \mathrm{d}u + \frac{1}{\beta_1} \sum_{r \in \mathcal{R}}\sum_{s \in \mathcal{S}} q_{rs}\left(\ln q_{rs} - 1 - \beta_2\rho_s - \beta_{0,s}\right)\label{obj:lower} \\
& & & \text{subject to} \nonumber & & \\
& & & & & \bm{v} = \sum_{r \in \mathcal{R}, s \in \mathcal{S}}\bm{\hat{v}}_{rs}+ \sum_{s \in \mathcal{S}, k \in \mathcal{K}}\bm{\check{v}}_{sk} , \; \label{cons:v}\\
& & & (\hat{\bm{\lambda}}_{rs}) & & A\bm{\hat{v}}_{rs} = q_{rs}E_{rs} , \; \forall r, s\label{cons:v_q}\\
& & & (\check{\bm{\lambda}}_{sk}) & & A\bm{\check{v}}_{sk} = \bar{d}_{sk}E_{sk} , \; \forall s, k\label{cons:v_d}\\
& & & ({\gamma}_{r})& &  \sum_{s \in \mathcal{S}}{q_{rs}} = Q_r, \; \forall r\label{cons:q_Q}\\
& & & ({\eta}_{sk})& &  \bar{d}_{sk} = \delta_{sk} \bar{D}_s,\; \forall s, k\label{cons:d_D}
\end{align}
\label{model:lower}
\end{subequations}

where: 
\begin{description}[leftmargin=!,labelwidth=\widthof{12345}]
\itemsep-0.7em
\item[$\bm{\hat{v}}_{rs}, \bm{\check{v}}_{sk}$] : link flow vector (each row corresponds to a link, $a$) grouped by OD for drivers and passenger vehicles, respectively;
\item[$\bm{q}$] : drivers' relocation flow vector (each row corresponds to an OD, $rs$);
\item[$v_a$] : aggregate traffic flow on link $a$; 
\item[$t_a(\cdot)$] : travel time function of link $a$; 
\item[$A$] : node-link incidence matrix of network, with $1$ at starting node and $-1$ at ending node;
\item[$E_{rs}$] : O-D incidence vector of O-D pair $rs$ with $1$ at origin $r$, $-1$ at destination $s$. If $r=s$, $E_{rs} = \bm{0}$;
\item[$\bar{d}_{sk}$] : passenger vehicle traffic flow from $s$ to $k$;
\item[$\delta_{sk}$] : ratio of traffic going from $s$ to $k$ to total travel originating from $s$ (model input); 
\item[$Q_{r}$] : drivers initial availability at $r$ (model input); 
\item[$\bar{D}_{s}$] : number of vehicles departing from $s$, including conventional passenger vehicles and ride-sourcing vehicles (model input); 
\item[$\check{\lambda}, \hat{\lambda}, \gamma, \eta$] : dual variables of the corresponding constraints. 
\end{description}

Constraint (\ref{cons:v}) calculates the aggregate link flow $v_a$ from the link flow associated with OD pairs $rs$ and $sk$. Constraints (\ref{cons:v_q}, \ref{cons:v_d}) are link-based formulations of flow conservation at each node for each OD pair. Constraint (\ref{cons:q_Q}) assumes the total drivers availability $Q_r$ at each location $r$ is inelastic, even though the relocation demand of drivers to each location $s$ is elastic and depends on the factors described in (\ref{eq:utility}). Constraint (\ref{cons:d_D}) calculates the distribution of riders and other conventional passenger vehicles OD demand.

In the objective function (\ref{obj:lower}), the first term corresponds to the total user cost as modeled in a conventional static traffic equilibrium model \citep{beckmann1956studies}, the second term involving $q^{rs}\ln q^{rs}$ corresponds to the entropy of trip distribution, and the remaining terms correspond to the utility measure of the drivers.  This objective function does not have a physical interpretation, as pointed out by \cite{Sheffi_85}, but it guarantees the first Wardrop principle \citep{Wardrop_52} and the multinomial logit location choice assumption being satisfied, as formally stated in {\bf Lemma \ref{lem:CDA}}.

\begin{lemma}{(combined distribution and assignment for ride-sourcing system)}\label{lem:CDA}
The optimal solutions $(\bm{\hat{v}^\ast},\bm{\check{v}^{\ast}}, \bm{q^\ast} )$ of  problem (\ref{model:lower}) are the equilibrium solutions for Wardrop user equilibrium and the relocation choice of drivers with multinomial logit model (\ref{eq:utility}) given spatial prices $\bm{\rho}$.
\end{lemma}  

\state Proof. See \ref{app:pfs}. 
\eop

\subsubsection{TNC Behaviors}

The pricing decisions of a TNC will influence the equilibrium solutions of (\ref{model:lower}), which forms a mathematical programming with equilibrium constraints (MPEC) problem. The objective of a TNC could be balancing supply and demand, maximizing matches, profits, market shares, services quality, minimizing total relocation distance, etc. In this paper, we assume that a TNC would like to minimize the imbalance between supply and demand \footnote{Balancing supply and demand is the goal of dyanmaic pricing claimed by all  TNC companies. See, for example, https://www.ridester.com/training/lessons/surge-primetime-boost-primezones/ (visited on Sep. 26, 2019).}, which has the advantage to improve driver/rider satisfactory and ultimately leads to a long-term sustainable development. Comparison between different objectives (minimizing imbalance versus maximizing short-term profits) will be presented in Section \ref{sec:nume}. The TNC decision making can be modeled as the upper-level problem, as formulated in (\ref{model:upper}).

\begin{subequations}
\begin{align}
& \underset{\boldsymbol{\rho} \in \reals^{S}}{\text{minimize}}
& &\sum_{s \in \mathcal{S}} m_s \label{obj:upper}\\
& \text{subject to}  & & \nonumber \\
& && m_s = |\sum_{r \in \mathcal{R}} q_{rs}- ({D}_s - b_s\rho_s)|, \ \forall s \label{cons:matches}\\
& && (\ref{model:lower})
\end{align}
\label{model:upper}
\end{subequations}
where: 
\begin{description}[leftmargin=!,labelwidth=\widthof{12345}]
\itemsep-0.7em
\item[$m^{s}$] : demand-supply imbalance at $s$.
\end{description}

Obejective (\ref{obj:upper}) minimizes the total demand-supply imbalance over the study area. Constraint (\ref{cons:matches}) calculates the imbalance as the difference between the number of riders requesting services ${D}_s - b_s \rho_s$ and the number of drivers arriving at this location $\sum_{r \in \mathcal{R}}q_{rs}$. While constraint (\ref{cons:matches}) can be reformulated as linear constraints, the bi-level structure still makes the whole problem challenging to solve to global optimal \citep{yang1998models, gao2005solution, sinha2017review}. In Section \ref{sec:solu_meth}, we propose an equivalent single-level convex reformulation of (\ref{model:upper}), which can be efficiently solved by commercial nonlinear solver (e.g. IPOPT \citep{wachter2009short}) or by classic traffic assignment algorithms (e.g. Frank-Wolfe algorithm \citep{Sheffi_85} or Evans' procedure \citep{evans1976derivation}). To investigate the transportation impacts, the spatial pricing models presented in this paper mainly focus on peak hours, so that we do not explicitly address the unmatched, if any, riders and drivers in the following time steps. The unserved customers could keep requesting ride-sourcing services, or switch to other transportation modes, such as driving. The idle ride-sourcing vehicles will decide whether to stay, relocate, or leave the platform. Their decisions will determined the drivers initial availability and riders demand function for the next time step.

\subsection{Mathematical Modeling: Considering Matching Time}\label{sec:modeling_matching}

In this section, we explicitly model the impacts of matching time. Similar to \cite{zha2016economic}, we assume the matching rate $m$ has constant elasticity with respect to vacant drivers $N_D$ and waiting riders $N_R$, which leads to a Cobb-Douglas matching function $m = \alpha_0(N_D)^{\alpha_1}(N_R)^{\alpha_2}$. In contrast to \citep{zha2016economic}, we do not assume a stationary state, i.e. we assume driver arrival rate, rider requesting rate, and matching rate are close but may not be equal during our studied periods. This relaxation could be helpful as we can see in Section \ref{sec:nume} that dynamic pricing does not guarantee supply-demand balance even in theory.

Denote the locational \footnote{We omit location index for simplicity.} average and total waiting time for drivers (D) and riders (R) as $t_D$, $T_D$, $t_R$, and $T_R$, respectively.  $\forall i \in \{D, R\}$, the average waiting time for $i$ can be calculated in (\ref{eq:waiting}), which is explained graphically in \ref{app:explain}.
\begin{subequations}
\begin{align}
& \Dot{N}_i(t) = f_i - m(t) = f_i - \alpha_0N_D^{\alpha_1}(t)N_R^{\alpha_2}(t), \label{eq:vacant}\\
& T_i = \int_{0}^{T}N_i(t) d{t}, \ \forall i \in \{D, R\} \label{eq:total_time}\\
& t_i = \frac{T_i}{f_iT} = \frac{\int_{0}^{T}N_i(t) d{t}}{f_iT},  \label{eq:average_time}
\end{align}
\label{eq:waiting}
\end{subequations}
where $f_{D}$ ($f_{R}$) is the ride-sourcing supply (demand) flow rate, and $T$ is the study period (e.g., peak hour). (\ref{eq:vacant}) is the dynamics of vacant drivers and waiting riders; (\ref{eq:total_time}) and  (\ref{eq:average_time}) calculate the total and average waiting time for drivers and riders over $T$. While solving (\ref{eq:waiting}) analytically is challenging, $t_i$ ($\forall i \in \{D, R\}$) can be solved numerically given $f_D, f_R$ and boundary conditions $N_D(0), N_R(0)$. We then apply linear regression techniques to estimate the relationship between $t_i$ and $f_i$, $\forall i \in \{D, R\}$, with a function form shown in (\ref{eq:waiting_time_flow}). When $\alpha_0 = 0.1, \alpha_1 = \alpha_2 = 0.6$, we have $a_0 = 6.29, a_1 = 2.24, a_2 = -2.40$ as least square estimators. $|a_2| > |a_1|$ implies increasing return of scale for matching rate. The regression summary output (Figure \ref{fig:regression_output}) is presented in \ref{app:explain}. 
\begin{equation}
t_{i} = a_0 f_{i}^{a_1}f_{-i}^{a_2}\quad \forall i \in \{D, R\}\label{eq:waiting_time_flow}
\end{equation}

\begin{figure*}[htbp]
	\centering
	\makebox[1.0\linewidth][c]{
		\begin{subfigure}[t]{.25\linewidth}
			\centering
			\caption{Original Network}
		\includegraphics[width=1\linewidth]{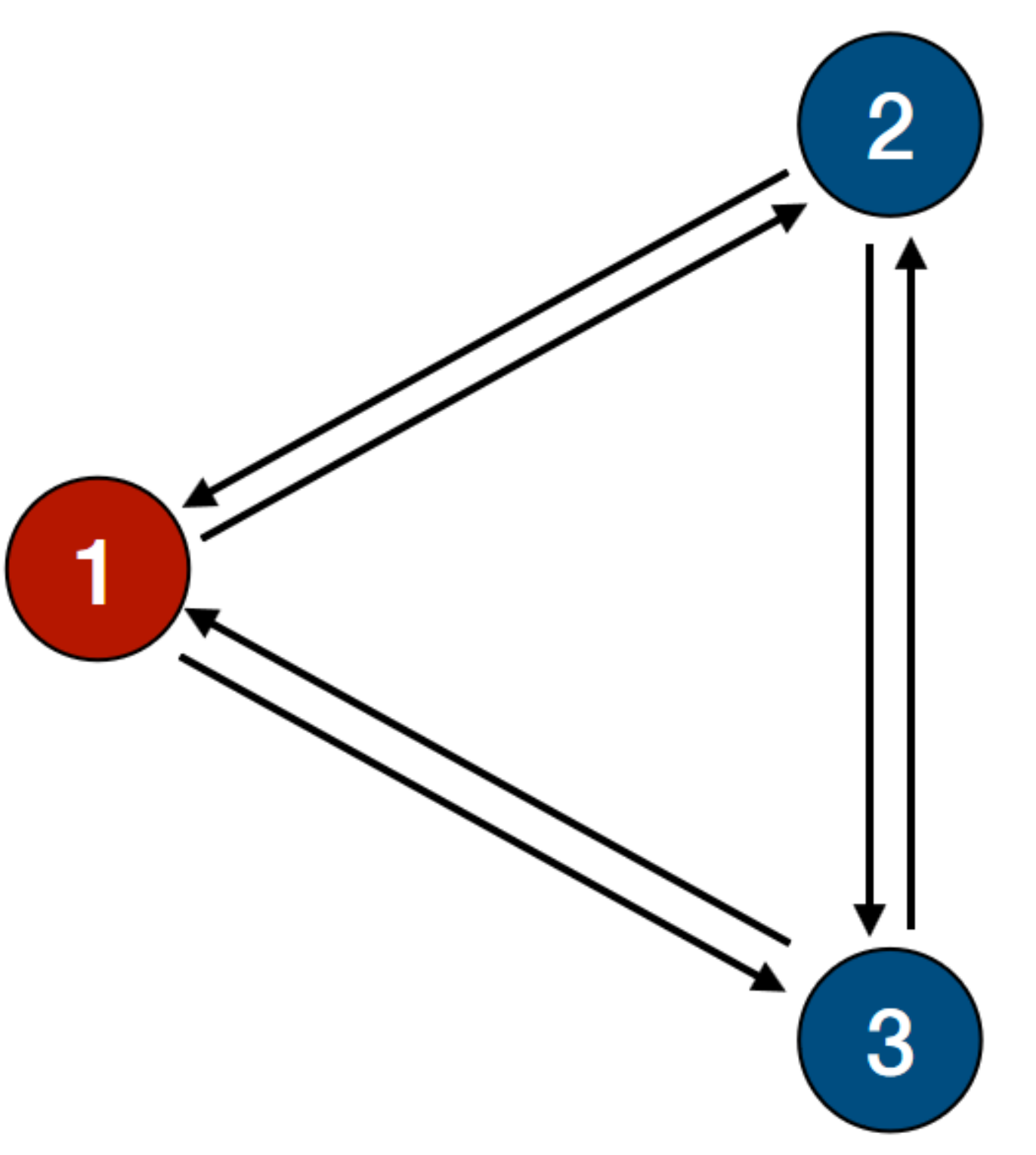}
			\label{fig:three_node}
		\end{subfigure}
		\begin{subfigure}[t]{.75\linewidth}
			\centering
			\caption{Augmented Network}
			\includegraphics[width=1\linewidth]{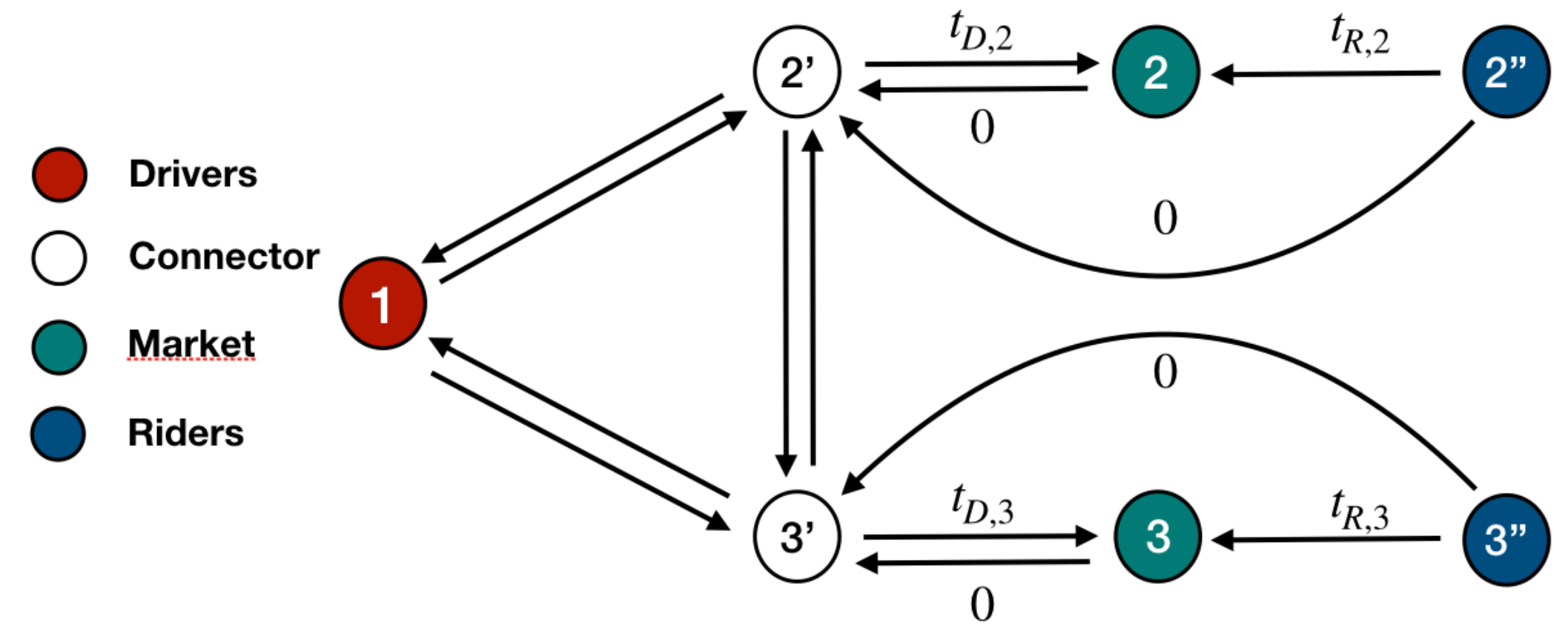} 
			\label{fig:seven_node}
	\end{subfigure}}
	\vspace{-2em}
	\caption{Illustration of Network Augmentation}
	\label{fig:aug_graph}
\end{figure*}

To incorporate matching time in a unified CDA framework, we define an augmented transportation network, $\bar{\mathcal{G}} = (\bar{\mathcal{N}}, \bar{\mathcal{A}})$. The augmenting procedures from  ${\mathcal{G}}$ to $\bar{\mathcal{G}}$ is illustrated in Figure \ref{fig:aug_graph} using a three-node transportation network, which has one driver node (node 1) and two rider nodes (node 2 and 3). For each rider node $s$ in Figure \ref{fig:three_node}, we created $s$, $s'$, and $s''$, representing nodes for markets, connectors, riders, respectively, with directed links connecting between them as shown in Figure \ref{fig:seven_node}. In augmented network $\bar{\mathcal{G}}$, link flow $f_{s's}$, $f_{s''s}$ and $f_{s''s'}$ represent the flow of driver supply, rider demand, and travelers choosing to drive at node $s$, while the link cost $t_{s's}$ and $t_{s''s}$ are the waiting time for drivers ($t_{D,s}$) and riders ($t_{R,s}$), respectively. If travelers choose to drive, there is no waiting time, so that link cost $t_{s''s'} = 0$. 


With an augmented network $\bar{\mathcal{G}}$, the utility function of drivers is identical with (\ref{eq:utility}). Note that $t_{rs}$ in (\ref{eq:utility}) will include not only the travel time from $r$ to $s'$, $t_{rs'}$, but also the matching time at $s$, $t_{s's}$. The drivers relocation and routing behaviors, therefore, can be modeled as (\ref{model:lower}) (replacing $\mathcal{G}$ with $\bar{\mathcal{G}})$, given riders demand $\bm{f}_{R}$ and spatial pricing $\bm{\rho}$.

For travelers at location $s$, the ride requesting demand function is not only depends on the prices $\rho_s$, but also depends on the waiting time, $t_{s''s}$. We assume that the choice of ride-sourcing over driving follows a binary logit model, with ride-sourcing utility function as (\ref{eq:utility_riders}),

\begin{equation}
\label{eq:utility_riders}
U_{s''s} = \beta_{0,s}' -\beta_1' t_{s''s} - \beta_2' \rho_{s}
\end{equation}
where $\beta_{0,s}', \beta_{1}', \beta_{2}'$ are the utility coefficients for ride-sourcing attractiveness, waiting time, and ride-sourcing prices. We assume driving and ride-sourcing will have the same travel time, therefore, travel time will be canceled out and does not need to be included in (\ref{eq:utility_riders}). Since the rider demand couples with waiting time (see (\ref{eq:waiting_time_flow}) and (\ref{eq:utility_riders})), the rider behaviors at each location $s$ can be formulated as a CDA model (\ref{model:lower_rider}) as well, where the first term corresponds to the conventional Wardrop user equilibrium, and the second term reflects the binomial logit choice between ride-sourcing and driving.

\begin{equation}
\underset{{f}_{R,s} \in \reals_+}{\text{minimize}}\quad \int_{0}^{f_{R,s}} t_{R,s}(f_{D,s},u) \mathrm{d}u + \frac{1}{\beta_1'} \left[f_{R,s}\left(\ln f_{R,s} - 1 + \beta_2'\rho_s - \beta_{0,s}'\right) + (D_s - f_{R,s})\left(\ln (D_s - f_{R,s}) - 1  \right)\right]
\label{model:lower_rider}
\end{equation}

Similar to formulation (\ref{model:upper}) for ignoring matching time case, the TNC decision making considering matching time can be modeled as bilevel problem, as shown in (\ref{model:upper_matching}), where objective (\ref{obj:upper}) minimizes total demand-supply imbalance (calculated by (\ref{cons:matches})), subject to drivers, riders, and all other travelers user equilibrium ((\ref{model:lower}, \text{replacing $\mathcal{G}$ with $\bar{\mathcal{G}}$}) and  (\ref{model:lower_rider})). 

\begin{subequations}
\begin{align}
& \underset{\boldsymbol{\rho} \in \reals^{S}}{\text{minimize}}
& &(\ref{obj:upper}) \\
& \text{subject to}  & & \nonumber \\
& && (\ref{cons:matches}) ,  (\ref{model:lower}, \text{replacing $\mathcal{G}$ with $\bar{\mathcal{G}}$}) ,  (\ref{model:lower_rider}) 
\end{align}
\label{model:upper_matching}
\end{subequations}

\section{Solution Approach}\label{sec:solu_meth}




\subsection{Solution Approach: Ignoring Matching Time} \label{sec:solu_approach_ignore}
We start by showing that when problem (\ref{model:upper}) is solved to optimal, the ride-sourcing market will be balanced at each location. In addition, optimal solutions $\bm{\rho}^*$ is unique. This is more formally stated in Lemma \ref{lem:balance_matches}.

\begin{lemma}{(balancing ride supply and demand)} \label{lem:balance_matches}
$\bm{\rho}^*$ optimizes problem (\ref{model:upper}) if and only if ride supply and demand are balanced at each location given $\bm{\rho}^*$. In addition, $\bm{\rho}^*$ is unique.
\end{lemma}  

\state Proof. See \ref{app:pfs}. 
\eop

Lemma \ref{lem:balance_matches} says that problem (\ref{model:upper}) is equivalent to the problem of finding the unique $\bm{\rho}^{\ast}$ that balance supply and demand at each location. Notice that $\bm{\rho}$ can be negative under extreme scenarios, when too many drivers exist in the system and TNC may provide incentives for riders to use ride-sourcing services. However, this will not happen in reality given reasonable driver supply $Q_r$, since drivers will dynamically respond to extreme low prices\citep{bimpikis2019spatial}. Leveraging Lemma \ref{lem:balance_matches}, the optimal spatial prices can be obtained by solving a single-level convex optimization problem, which will be illustrated in the remaining of this section.

%

Problem (\ref{model:upper}) has three types of agents: drivers, riders, and TNC. Because of Lemma \ref{lem:balance_matches}, the decision making of TNC can be reformulated as finding $\bm{\rho}^*$ such that 

\begin{equation}
\label{cons:market_clearing}
    \sum_{r \in \mathcal{R}} q_{rs}^{\ast}(\bm{\rho}^*) = D_s - b_s\rho_{s}^{\ast}, \forall s \in \mathcal{S}
\end{equation}
where $q_{rs}^{\ast}(\bm{\rho}^*)$ is the optimal solutions of (\ref{model:lower}) given $\bm{\rho}^*$. This problem falls within the framework of multi-agent optimization problem with equilibrium constraints (MOPEC) \citep{ferris2013mopec}, where drivers make relocation decisions; riders decide whether to request ride-sourcing services; all travelers make route choices decisions; and the system is subject to equilibrium conditions (\ref{cons:market_clearing}). MOPEC has been applied in transportation, energy, economics domains, such as charging infrastructure planning \citep{Guo_et_al_16}, hydro-thermal electricity systems optimization \citep{philpott2016equilibrium}, renewable energy supply chain planning \citep{guo2017stochastic}, and computational Walrasian equilibrium \citep{deride2019solving}.

Define two types of dummy agents in our problem, $Driver$ and $Rider$. Given $\bm{\rho}$, $Driver$ solves problem (\ref{model:lower}), which determines the drivers relocation decision from $r$ to $s$, $q_{rs}(\bm{\rho}), \forall r, s$; and $Rider$ solves problem (\ref{model:demand}), which determines the riders demand $d_s(\bm{\rho})$ at each location $s$. Note that problem (\ref{model:demand}) is a constructed unconstrained quadratic programming problem, which has a closed-form solution $d_s^{\ast}(\bm{\rho}) = D_s - b_s \rho_s, \forall s \in \mathcal{S}$, identical to the demand function (\ref{eq:demand}) we assumed.  

\begin{equation}
\label{model:demand}
\underset{\bm{d} \in \reals^{S}}{\text{maximize}} \ \ \sum_{s \in \mathcal{S}} D_s d_s - b_s\rho_s d_s - \frac{1}{2}(d_s)^2
\end{equation}

The equilibrium constraints (\ref{cons:market_clearing}) can be reformulated as (\ref{model:equilibrium_constraints}). 

\begin{equation}
\label{model:equilibrium_constraints}
\sum_{r \in \mathcal{R}} q_{rs}^{\ast}(\bm{\rho}) = d_s^{\ast}(\bm{\rho}), \forall s \in \mathcal{S}
\end{equation}

Notice that both the objective functions of (\ref{model:lower}) and (\ref{model:demand}) consist of terms regarding prices times quantities, i.e., $\sum_{r \in \mathcal{R}}\rho_sq_{rs}$ and $\rho_sd_s$. In addition, notice that the equilibrium constraints (\ref{model:equilibrium_constraints}) are quantity balance of $\sum_{r \in \mathcal{R}}q_{rs}$ and $d_s$. Inspired by Lagrange multiplier method, we construct a single-level convex optimization problem (\ref{model:combine}), which will yield the same spatial pricing decisions as bilevel formulation, as stated in Theorem \ref{thm:combine}.

\begin{theorem}{(single-level convex reformulation)}\label{thm:combine}
$\bm{\rho}$ solves bi-level problem (\ref{model:upper}) if and only if $\bm{\rho}$ solve single-level problem (\ref{model:combine}).
\end{theorem}

\state Proof. See \ref{app:pfs}. 

\begin{subequations}
\begin{align}
& & & \underset{\bm{\hat{v}},\bm{\check{v}}, \bm{q}, \bm{d} \in \reals_+}{\text{minimize}}
& & \frac{\beta_1}{\beta_2}\sum_{a \in \mathcal{A}} \int_{0}^{v_a} t_a(v_a) \mathrm{d}u + \frac{1}{\beta_2} \sum_{r \in \mathcal{R}}\sum_{s \in \mathcal{S}} q_{rs}\left(\ln q_{rs} - 1 - \beta_{0,s}\right) + \sum_{s \in \mathcal{S}}\frac{1}{b_s}(\frac{d_s^2}{2} - D_sd_s)\label{obj:combine} \\
& & & \text{subject to} \nonumber & & \\
& & & & & (\ref{cons:v} \sim \ref{cons:d_D}) \nonumber\\
& & & ({\rho}_{s})& &  \sum_{r \in \mathcal{R}} q_{rs} = d_s,\; \forall s\label{cons:q_d}
\end{align}
\label{model:combine}
\end{subequations}

The objective function (\ref{obj:combine}) is a combination of scaled objective functions (\ref{obj:lower}) and (\ref{model:demand}) without $\sum_{r \in \mathcal{R}}\rho_sq_{rs}$ and $\rho_sd_s$ terms. The purpose of scaling (\ref{obj:lower}) and (\ref{model:demand}) is to convert the units of both objective functions into \$ so that they are addable. In addition, the optimal dual variables $\rho_s^{*}$ of constraints (\ref{cons:q_d}) will unit in \$ and balance supply and demand at each location. Therefore, $\bm{\rho}^{*}$ can be interpreted as spatial prices to solve problem (\ref{model:upper}) based on Lemma \ref{lem:balance_matches}.

Problem (\ref{model:combine}) is a non-linear convex problem, with convex objective function and linear constraints. So it can be solved to global optimal by commercial nonlinear solver, such as IPOPT, relatively efficiently. For extreme large-scale problems, this reformulation opens up opportunities to apply classic and efficient transportation network solution algorithm, such as Frank-Wolfe algorithm \citep{Sheffi_85} and Evans' procedure \citep{evans1976derivation}. In this paper, we focus on small/medium transportation networks to draw policy insights, and leave the development of advanced algorithms for large network in the future.

\subsection{Solution Approach: Considering Matching Time}

Directly solving TNC upper-level problem (\ref{model:upper_matching}) may be challenging due to highly non-convexity. We restate (\ref{model:upper_matching}) as a problem of finding a maxinf-point for an appropriate bivariate function (bi-function) $W(\bm{\rho},\bm{\phi})$. This conversion is inspired by a recent theoretical development on variational convergence of bifunction by \cite{jw_09}, which offers the flexibility of constructing a bifunction and allows one to choose a sequence of bifunctions with desired properties (including convexity and continuity) to approximate the original bifunction.  A similar approach has been successfully implemented in \citep{Guo_et_al_16, deride2019solving}. In this paper, we highlight the key steps of the algorithm design.

We introduce the \emph{Walrasian} 
function associated with this equilibrium problem, defined as

\[W(\bm{\rho},\bm{\varphi})=-\sum_{s\in \mathcal{S}} \varphi_s{\rm ES}_s(\bm{\rho}),\quad {\rm on }\, \reals^S\times\Delta_S,\]

\noindent where $\Delta_S$ corresponds to the $S$-dimensional unit simplex. $\rm{ES_s}$ is the square of demand-supply imbalance, defined in (\ref{eq:excess_supply}).

\begin{equation}
{\rm ES}_s(\bm{\rho}) = \left[\sum_{r \in \mathcal{R}}q_{rs}^{\ast}(\bm{\rho}) - f_{R,s}^{\ast}(\bm{\rho})\right ]^2, \; \forall s \in \mathcal{S}
\label{eq:excess_supply}
\end{equation}
where $q_{rs}^{\ast}$ and $f_{R,s}^{\ast}$ are the optimal solutions (as a function of $\bm{\rho}$) of problem (\ref{model:lower}, replacing $\mathcal{G}$ with $\bar{\mathcal{G}}$) and (\ref{model:lower_rider}), respectively. 

The relationship between maxinf point of $W(\bm{\rho},\bm{\varphi})$ and the demand-supply balancing (equilibrium) prices is stated in Lemma \ref{lem:mieq}.

\begin{lemma}{(equilibrium prices and maxinf-points)}\label{lem:mieq}
For $\bm{\rho} \in\reals^S$, $\bm{\rho}$ is a maxinf-point of the Walrasian function $W$ such that $W(\bm{\rho}, \cdot) \geq 0$,
on $\Delta_S$ if and only if $\bm{\rho}$ is an equilibrium point. 
\end{lemma}

\state Proof. See \ref{app:pfs}. \eop

Since $\inf_{\bm{\varphi}} W(\bm{\rho},\bm{\varphi})$ is lack of concavity in $\bm{\rho}$, we apply an \emph{augmented Lagrangian}
for this non-concave formulation. Given sequences of nonegative, nondecreasing scalars $\{r^\nu\},\,\{M^\nu\}$, one can define the sequence of \emph{augmented Walrasian}
functions for this problem as
\begin{equation}
    W^\nu(\bm{\rho},\bm{\varphi})=\inf_{\bm{z}}\left\{\left.W(\bm{\rho},\bm{z})+\frac{1}{2r^\nu}|\bm{z}-\bm{\varphi}|^2\right|\bm{z}\in\Delta_S\right\},{\rm on }\, [\bm{\rho}^{\nu -1}-M^\nu,\bm{\rho}^{\nu-1}+M^\nu]\times\Delta_S
\end{equation}

The idea of this procedure is to approximate the problem of finding maxinf-points of the original Walrasian function $W$, by computation of approximate maxinf-points $\bm{\rho}_\varepsilon$ (see Definition \ref{def:approx_maxinf}) given by a sequence
of augmented Walrasians $W^\nu$, which are easier to solve. The convergence theorem of the proposed approximation scheme is stated in Theorem \ref{thm:algo_lopconv}.

\begin{defn}{(approximating equilibrium point)}
$\bm{\rho}_\varepsilon$ is an 
$\varepsilon$-approximating equilibrium point, denoted as $\varepsilon-\nargmaxinf W$, for $\varepsilon\geq 0$, if the following inequality holds:
$\left|\ninf W(\bm{\rho}_\varepsilon,\cdot)- \nsupinf W\right|\leq \varepsilon.$
\label{def:approx_maxinf}
\end{defn}

\begin{theorem}{(convergence of approximating maxinf-points)}\label{thm:algo_lopconv}
Suppose that $\nsupinf W$ is finite. Consider non-negative sequences $\{r^\nu\}$, $\{M^\nu\}$, and $\{\varepsilon^\nu\}$
such that $r^\nu \upto \infty$, $M^\nu\upto \infty$, $\varepsilon^\nu\downto 0$. Let $\{W^\nu\}$
be a family of augmented Walrasian functions associated with each parameters $r^\nu$ and $M^\nu$.
Let $\bm{\rho}^\nu\in\varepsilon^\nu-\nargmaxinf W^\nu$,
and $\bm{\rho}^*$ be any cluster point of $\{\bm{\rho}^\nu\}$. Then $\bm{\rho}^*\in\nargmaxinf W$.
\end{theorem}

\state Proof. See \citep{Guo_et_al_16}. \eop





Following Theorem \ref{thm:algo_lopconv}, we propose the computational algorithm, as summarized in Algorithm \ref{alg:iter}, to achieve a sequence of approximating maxinf-points.

\begin{algorithm}[H]
	\KwResult{$\boldsymbol{\rho}^{\nu}$}
	initialization: $\nu = 0$, $\varepsilon^0$, $M^0$, $r^0$, ${\rm gap}^0$, $\boldsymbol{\rho}^0$, $\epsilon$, $0<c_1<1$, $c_2>1$\;
	
	\While{${\rm gap}^{\nu} \geq \epsilon$}{
	
	Phase I: solve the minimization problem $\bm{\varphi}^{\nu+1}\in\nargmin W^{\nu+1}(\bm{\rho}^\nu,\cdot)$\;
	
    Phase II: solve the maximization problem $\bm{\rho}^{\nu+1}\in\nargmax W^{\nu+1}(\cdot,\bm{\varphi}^{\nu+1})$\;
    
    evaluate: ${\rm ES}_s(\bm{\rho}^{\nu+1}), s \in \mathcal{S}$ \;
    
    let: $gap^{\nu+1} = \max\{{\rm ES}_s(\bm{\rho}^{\nu+1}), s \in \mathcal{S}\}$\;
    
	let: $\varepsilon^{\nu+1} = c_1\varepsilon^{\nu}$, $M^{\nu+1} = c_2M^{\nu}$, $r^{\nu+1} = c_2 r^{\nu}$\;
	
	let: $\nu := \nu + 1$
	}
	\caption{Approximating Maxinf-Point Algorithm}
	\label{alg:iter}
\end{algorithm}

Phase I consists of the minimization of a quadratic objective function over the $S$-dimensional simplex. This can be solved using Cplex/Gurobi solver. Phase II is done without considering first order information (due to a lack of concavity) and relying on BOBYQA algorithm \citep{Pow09:bobyqa}, which performs a sequentially local quadratic fit of the objective functions, over box constraints, and solves it using a trust-region method. As $r^\nu\upto \infty$, and $M^\nu\upto\infty$, and $\varepsilon^\nu\downto 0$, in virtue of Theorem \ref{thm:algo_lopconv}, $\bm{\rho}^\nu\to\bm{\rho}^*$, a maxinf-point of $W$.

\noindent{\bf Remarks}: 

1. The network augmenting procedures proposed in Section \ref{sec:modeling_matching} result in a transportation network with non-separable link-cost function (i.e., $t_{D}$ and $t_R$), and the Jacobian of link-cost function is asymmetric, for which there is no known mathematical program whose solutions are the equilibrium flow pattern \citep{Sheffi_85}. But, the traffic equilibrium solutions can be achieved by diagonalization method \citep{Sheffi_85}.  We note that since the Jocobian of the link-cost functions is not positive definite when $|a_1|\; \leq|a_2|$ and $f_{D,s} \approx f_{R,s}$, the convergence of traffic equilibrium is not guaranteed in theory, although the algorithm may still convergence when this condition is violated \citep{fisk1982solution}. To make our algorithm robust to non-convergence of traffic equilibrium, the algorithm we develop here does not require convergence for every iteration.

2. Notice that problem (\ref{model:upper_matching}) has at least $|\mathcal{S}|$ trivial equilibrium solutions, where only one location $s$ attracts all the drivers $\sum_{r \in \mathcal{R}} Q_r$. This is because when driver supply approaches 0 and price is finite at one location, rider demand will also approach 0 due to infinite waiting time, and vice versa. The insight from this is that ride-sourcing market may need a critical mass of drivers and riders to sustain. In this paper, we only focus on non-trivial equilibrium. 

\section{Numerical Examples} \label{sec:nume}
In this section, we test our models and solution approaches on a three-node system (Figure \ref{fig:three_nodes}) and Sioux Falls network (Figure \ref{fig:siou_fall}). We implement our model on Pyomo 5.6.6 \citep{hart2017pyomo}. We solve the problems using IPOPT 3.12.13 \citep{wachter2009short} for ignoring matching time case, and Cplex 12.8 \& BOBYQA (in NLopt 2.4.2) for considering matching time case, both with 0.1\% optimality gap. All the numerical experiments presented in this section were run on a 3.5 GHz Intel Core i5 processor with 8 GB of RAM memory, under Mac OS X operating system. 

\subsection{Three-node Network}
\begin{figure}[htbp]
\begin{center}
    \includegraphics[width=0.5\textwidth]{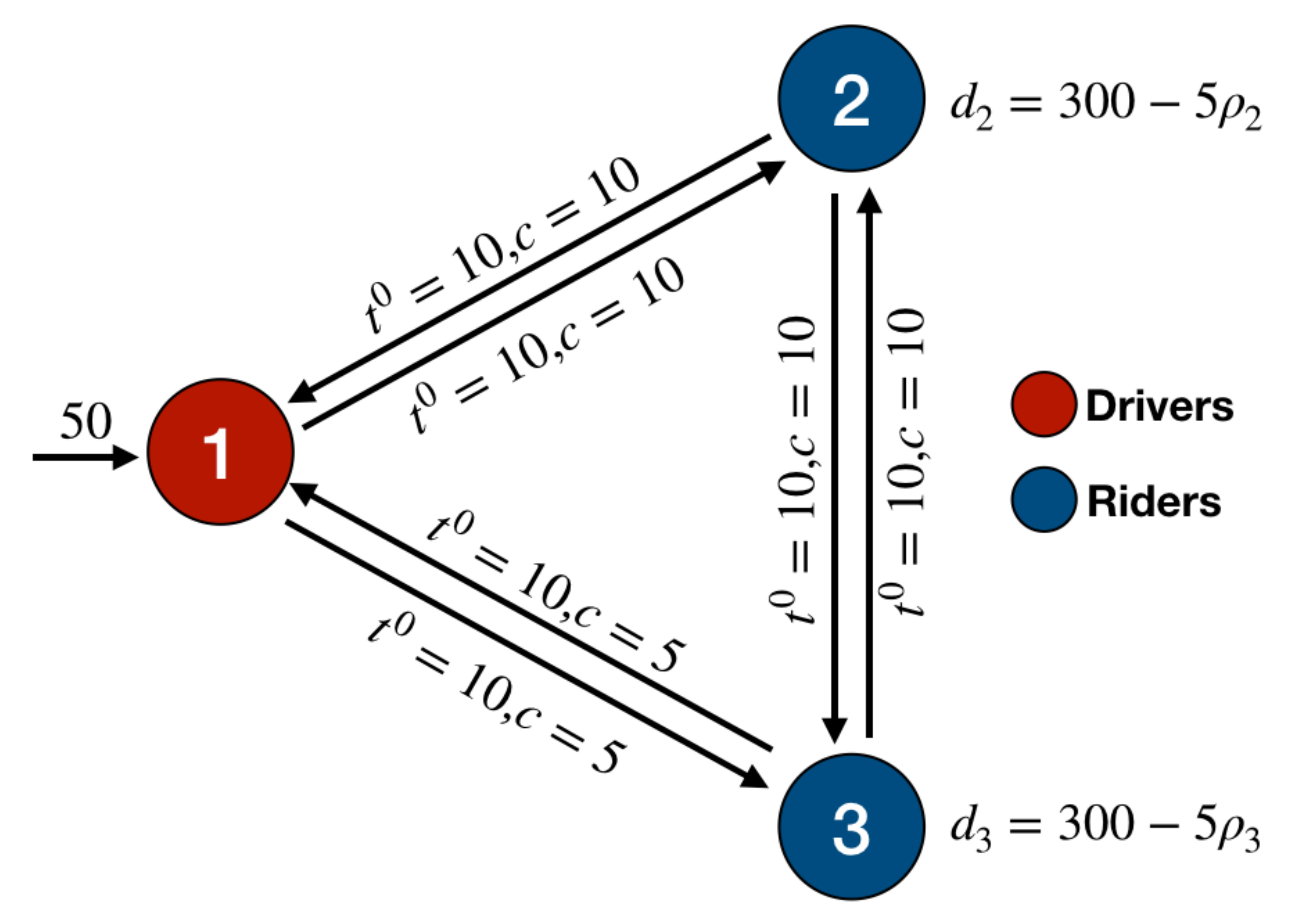}
\caption{Three Nodes Test Network}
\label{fig:three_nodes}
\end{center}
\end{figure}

The three-node test instance (Figure \ref{fig:three_nodes}) has one driver node (node 1) and two demand nodes (node 2 and 3). 50 drivers decide where to pick up riders at node 1. $\beta_0 = 0, \beta_1 = 1, \beta_2 = 0.6$. The link travel time has a form of $t_a = t_a^0[1+0.15*(v_a/c_a)^2]$, where $t_a^0$ is the free flow travel time and $c_a$ is the link capacity parameter. The value of $t_a^0$ and $c_a$ for each link are shown in Figure \ref{fig:three_nodes}. Notice that this test instance is symmetric except that links between node 1 and 3 have lower capacity. Through out Section \ref{sec:nume}, background traffic is not considered. 

\subsubsection{Ignoring Matching Time}
We assume both node 2 and 3 have a demand function of $d_s = 300 - 5\rho_s$. Solving this problem using our single-level reformuation (\ref{model:combine}), we have $\rho_2 = 53.5, \rho_3 = 56.5$.  Notice that with these locational prices, the demand market at both node 2 and 3 are balanced. The difference between $\rho_2$ and $\rho_3$ is because of the link congestion between node 1 and 3. In other words, node 3 is more congested to travel to due to limited link capacity, so that to attract more drivers to pick up riders at node 3, TNC needs to offer a higher locational price. If TNC doesn't consider transportation congestion, the optimal locational prices are $\rho_2' = \rho_3' = 55.0$, which will lead to unbalance of supply and demand on both demand nodes in reality. These optimal solutions are validated by applying global solver SCIP \citep{achterberg2009scip} directly on mixed-integer nonlinear programming (MINLP) reformulation of bi-level problem (\ref{model:upper}). Reformulating bi-level problem into MINLP is standard, so we omit the detail procedures for conciseness. For three-node example, our reformulation (\ref{model:combine}) can be solved by IPOPT in 0.1 seconds, while MINLP formulation takes SCIP 5.6 seconds to solve. Although both SCIP and IPOPT solve the problem for global optimal, MINLP formulation is sensitive to the selection of big-M parameters, which is a challenging issue for mixed integer programming with switching constraints \citep{guo2016contingency}.

Because TNC aims to balance supply and demand, transportation congestion can be worse under spatial pricing. We compare between two scenarios, surge pricing and uniform pricing. If we change from uniform pricing to surge pricing, the total travel time increases slightly from 1333 to 1336. The increase of total travel time is because TNC gives incentive (surge pricing) to drivers traveling through congested area in order to balance demand. This result matches the empirical findings in \citep{nie2017can} that ride-sourcing may increase transportation congestion compared to taxi services. Regulations will be needed to prevent pricing selfishly by TNCs to improve overall transportation mobility, especially in congested area.

Next, we compare the impacts of TNC's objectives between maximizing profits and minimizing imbalance. The problem of maximizing profits is shown in (\ref{model:upper2}).

\begin{subequations}
\begin{align}
& \underset{\boldsymbol{\rho} \in \reals^{S}}{\text{maximize}}
& &\sum_{s \in \mathcal{S}} \rho_s n_s \label{obj:upper2}\\
& \text{subject to}  & & \nonumber \\
& && n_s = \min(\sum_{r \in \mathcal{R}} q^{rs}, {D}^s - b^s\rho_s), \ \forall s \label{cons:matches2}\\
& && (\ref{model:lower})
\end{align}
\label{model:upper2}
\end{subequations}
where: 
\begin{description}[leftmargin=!,labelwidth=\widthof{12345}]
\itemsep-0.7em
\item[$n_{s}$] : ride-sourcing matches at $s$.
\end{description}

\begin{figure*}[htbp]
	\centering
	\makebox[1.0\linewidth][c]{
		\begin{subfigure}[t]{.5\linewidth}
			\centering
			\caption{Sensitivity on Supply}
			\includegraphics[width=1\linewidth]{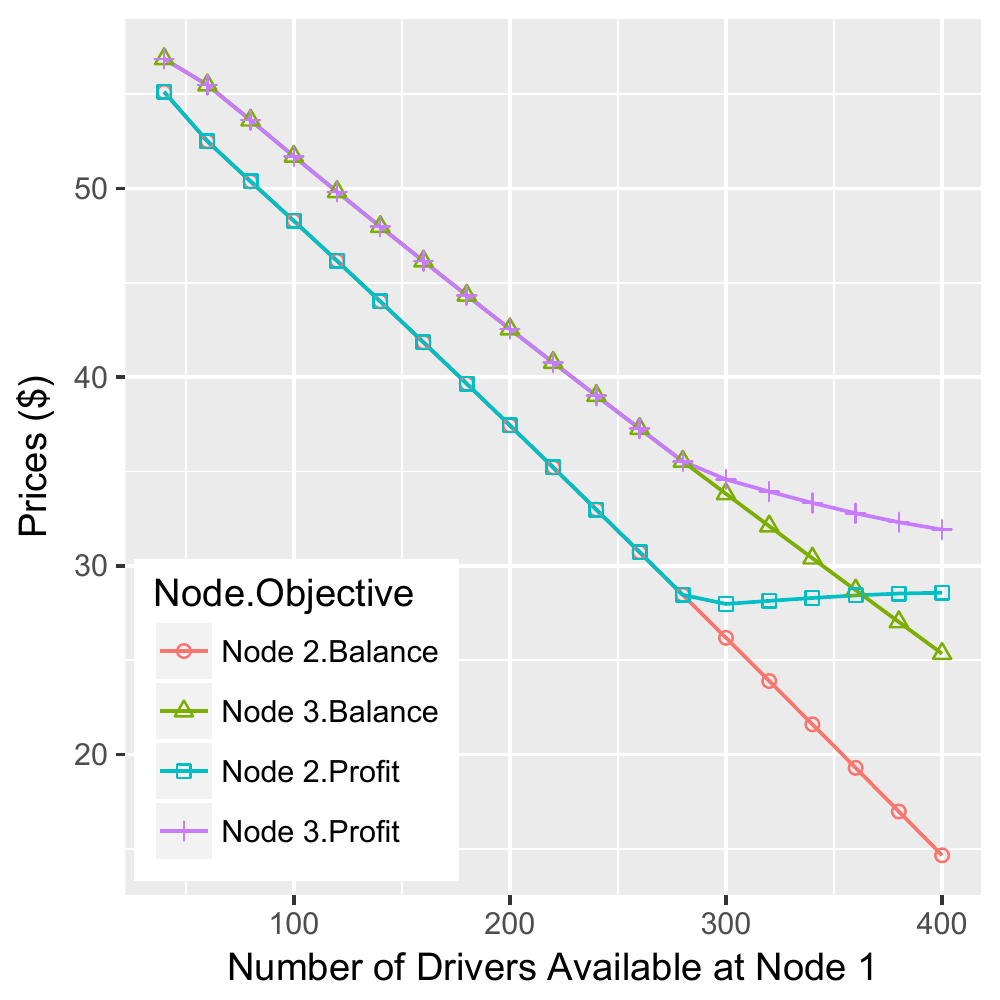}
			\label{fig:obj_supply}
		\end{subfigure}
		\begin{subfigure}[t]{.5\linewidth}
			\centering
			\caption{Sensitivity on Demand}
			\includegraphics[width=1\linewidth]{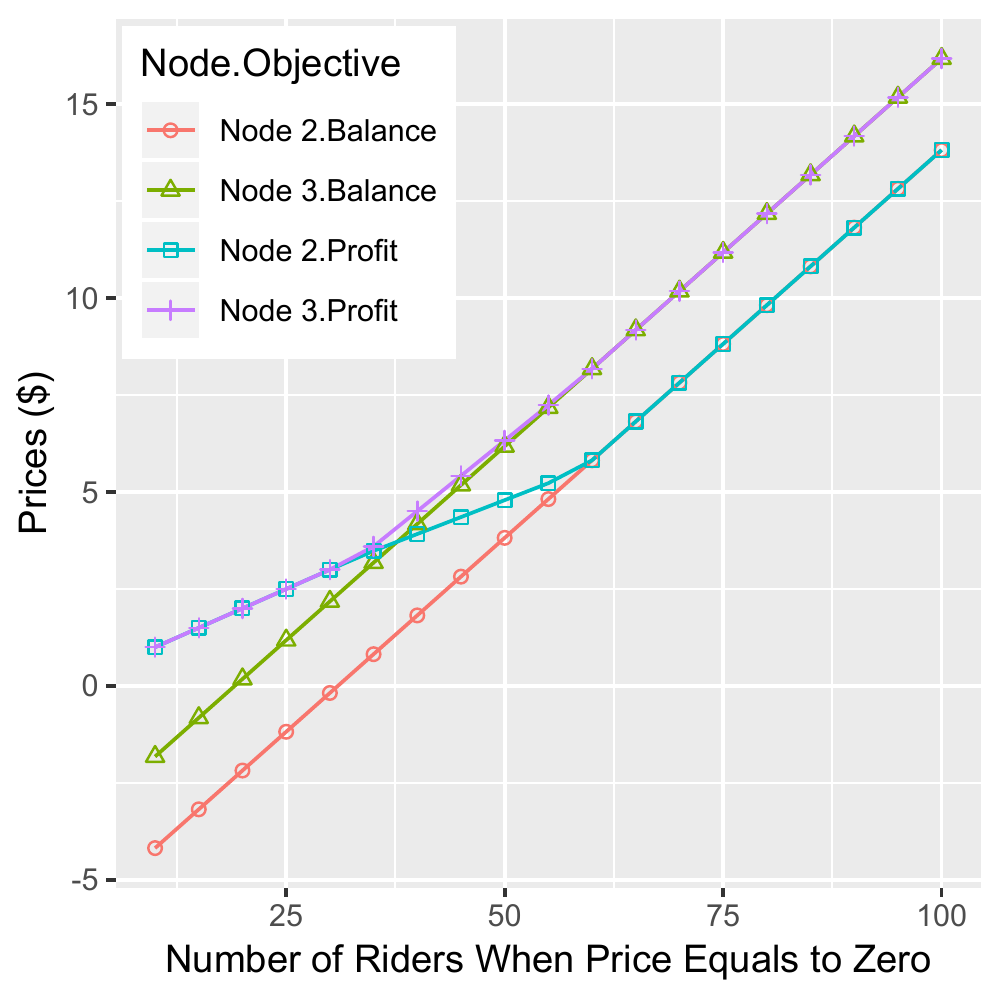} 
			\label{fig:obj_demand}
	\end{subfigure}}
	\vspace{-2em}
	\caption{Impacts of Supply and Demand on Optimal Surge Pricing}
	\label{fig:senstivity_price}
\end{figure*}

(\ref{obj:upper2}) maximizes total revenues. Since the majority of TNC operating costs go to sales, marketing, administrative, and R\&D \footnote{https://news.crunchbase.com/news/understanding-uber-loses-money/}, we assume that operational costs for TNC is independent of spatial pricing and drivers/riders flow, so the optimal pricing maximizes total revenues will also maximize total profits in the service territory. But we note that operating costs can be easily incorporated in (\ref{obj:upper2}) if operating costs data is available. (\ref{cons:matches2}) calculates the matches at each location $s$. For problem (\ref{model:upper2}), we are not able to directly make the single-level reformulation as in (\ref{model:combine}). So we solve the bilevel problem (\ref{model:upper2}) as MINLP using SCIP solver. The numerical results are shown in Figure \ref{fig:obj_supply} and Figure \ref{fig:obj_demand} , in which we conduct sensitivity analyses of the total number of drivers and interception of demand function on equilibrium prices with different TNC's objectives. When TNC aims to minimize imbalance for each location, the optimal prices are almost linearly decreasing (increasing) with increasing of supply (demand) for this toy test instance. The trend of optimal prices are slightly more complex when TNC aims to maximize profits. When supply is relatively high or demand is relatively low, TNC's profits are not restricted by drivers availability, so that TNC will set prices converge to monopoly prices, which are $D^s/2b^s, \forall s$. This means that when TNCs aim to maximize profits, the supply and demand may not necessarily balanced at each location. This finding is consistent with the results in \citep{bimpikis2019spatial} where endogenously determined spatial prices aiming to maximize profits may not always balance supply and demand. Notice that when demand is comparable to supply, prices minimizing imbalance of ride requests are identical with prices maximizing TNC profits, which indicates that maximizing balance of supply and demand leads to maximum profits for TNC.

\subsubsection{Considering Matching Time}
When matching time cannot be ignored for both drivers' and riders' decision making, we can solve problem (\ref{model:upper_matching}) using Algorithm \ref{alg:iter}. We assume $D_s = 300, \; \forall s$. The riders' utility coefficients $\beta_0' = 0, \beta_1' =1, \beta_2' = 0.6$. 

\begin{figure*}[htbp]
	\centering
		\begin{subfigure}[t]{.49\linewidth}
			\centering
			\caption{$\bm{\rho}^{0} \in [10, 20]$}
			\includegraphics[width=1\linewidth]{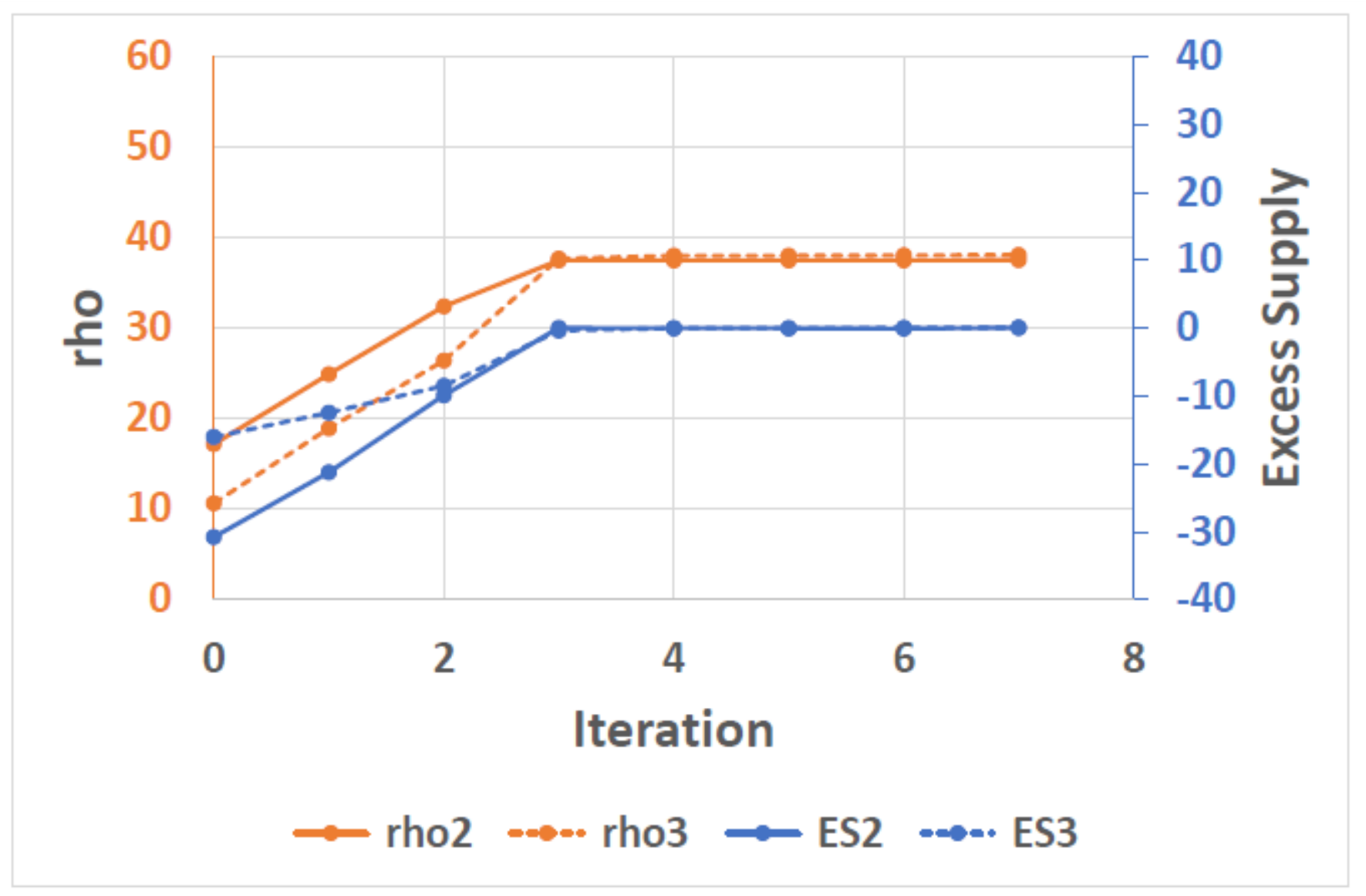}
			\label{fig:rho_1020}
		\end{subfigure}
		\begin{subfigure}[t]{.49\linewidth}
			\centering
			\caption{$\bm{\rho}^{0} \in [20, 30]$}
			\includegraphics[width=1\linewidth]{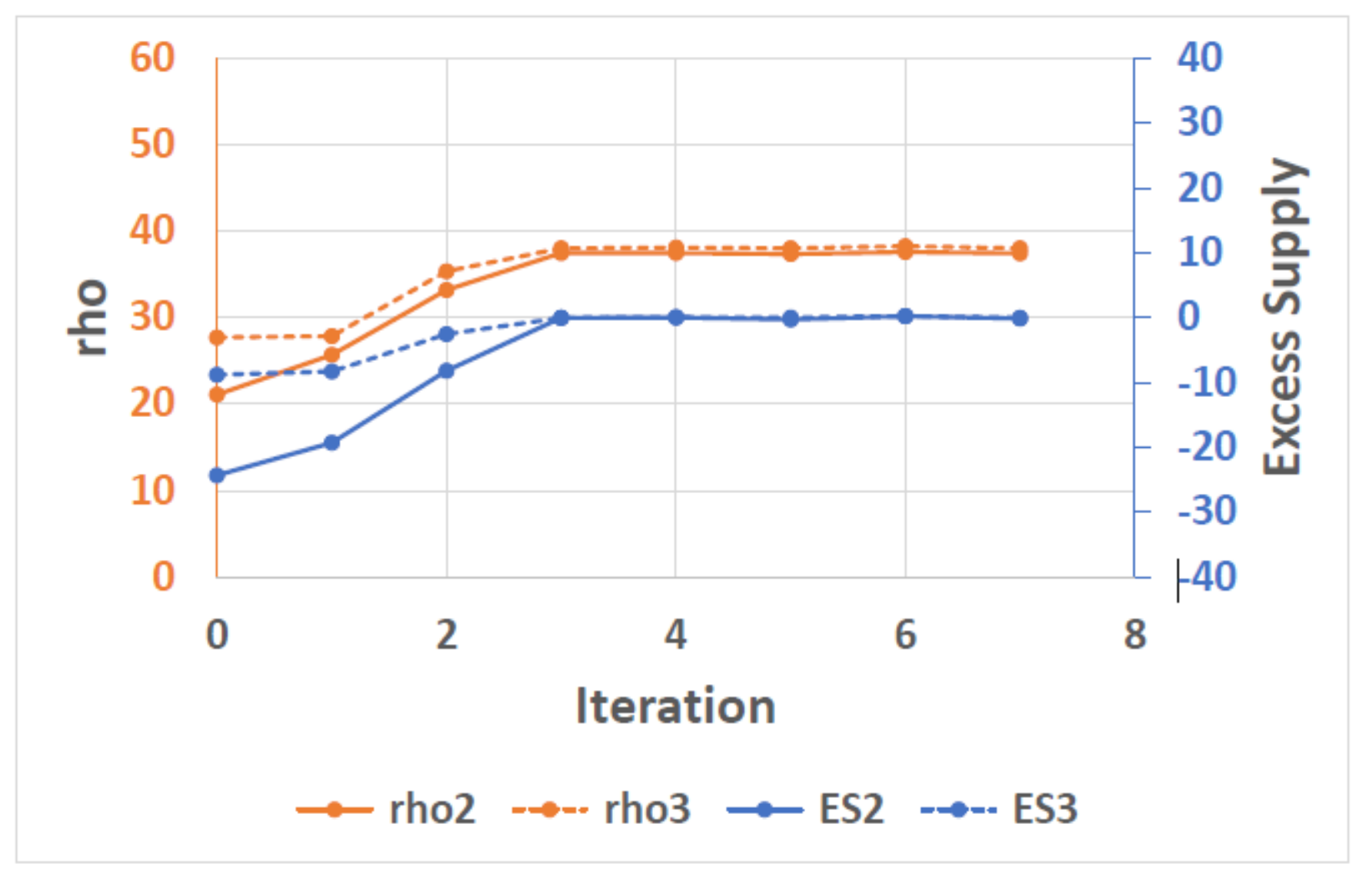} 
			\label{fig:rho_2030}
	    \end{subfigure}
	    \begin{subfigure}[t]{.49\linewidth}
			\centering
			\caption{$\bm{\rho}^{0} \in [30, 40]$}
			\includegraphics[width=1\linewidth]{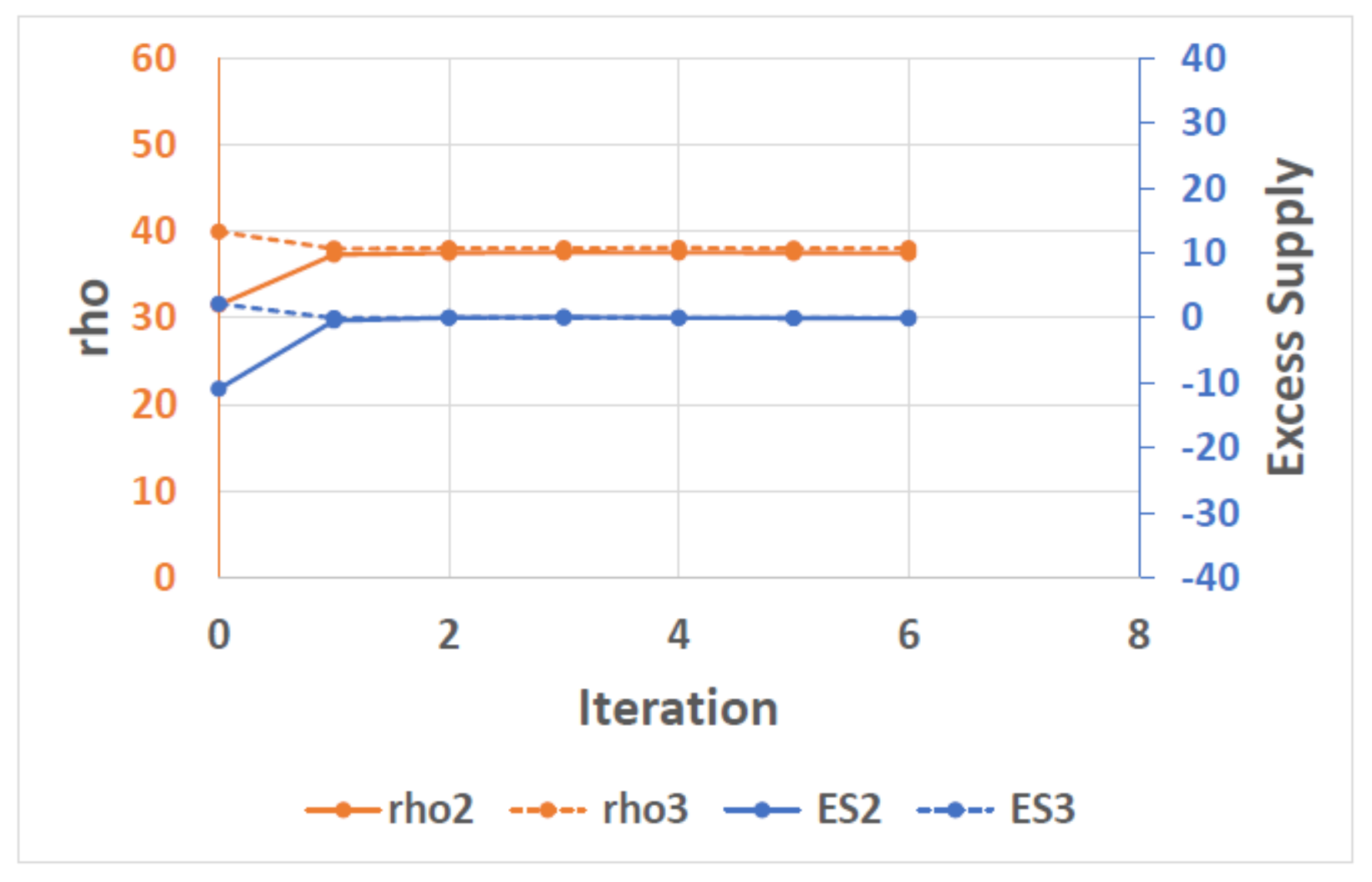}
			\label{fig:rho_3040}
		\end{subfigure}
		\begin{subfigure}[t]{.49\linewidth}
			\centering
			\caption{$\bm{\rho}^{0} \in [40, 50]$}
			\includegraphics[width=1\linewidth]{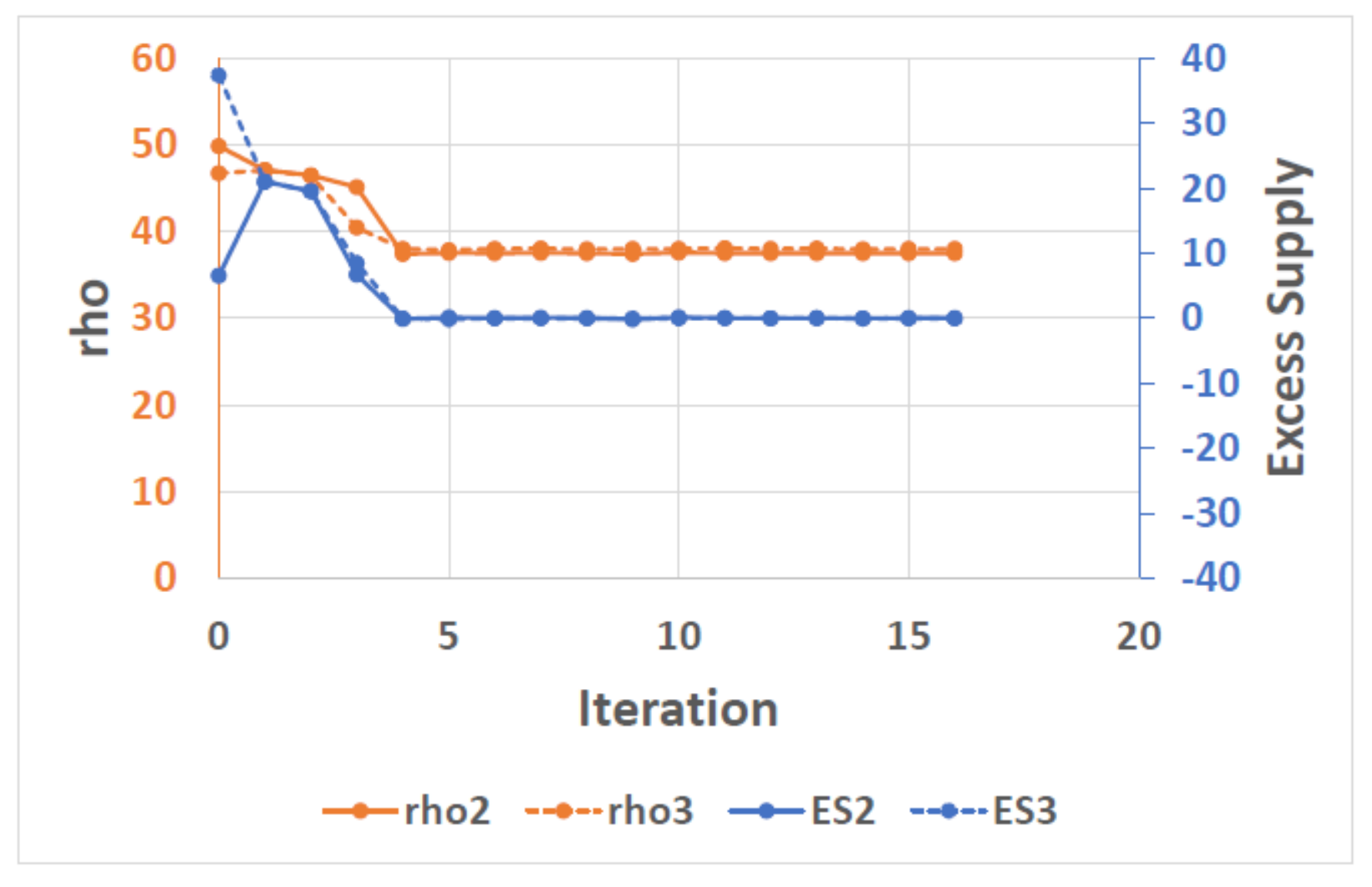} 
			\label{fig:rho_4050}
	\end{subfigure}
	\vspace{-2em}
	\caption{Convergence Patterns with Different Starting Points $\bm{\rho}^{0}$}
	\label{fig:three_node_conv}
\end{figure*}

With different initial locational prices, the convergence pattern of prices and ES are shown in Figure \ref{fig:three_node_conv}. The average computation time is 27.8s, and all of the experiments converged efficiently to the same equilibrium solutions, $\rho_2^* = 37.4, \rho_3^* = 38.0, f_{D,2}^*=f_{R,2}^* = 33.1, f_{D,3}^* = f_{R,3}^* = 16.9$. Notice that TNC offers a slightly higher price at location 3, which is more congested to travel to. This observation is consistent with the observation in ignoring matching time case. 

\begin{figure*}[htbp]
	\centering
	\makebox[1.0\linewidth][c]{
		\begin{subfigure}[t]{.333\linewidth}
			\centering
			\caption{Drivers Supply}
			\includegraphics[width=1\linewidth]{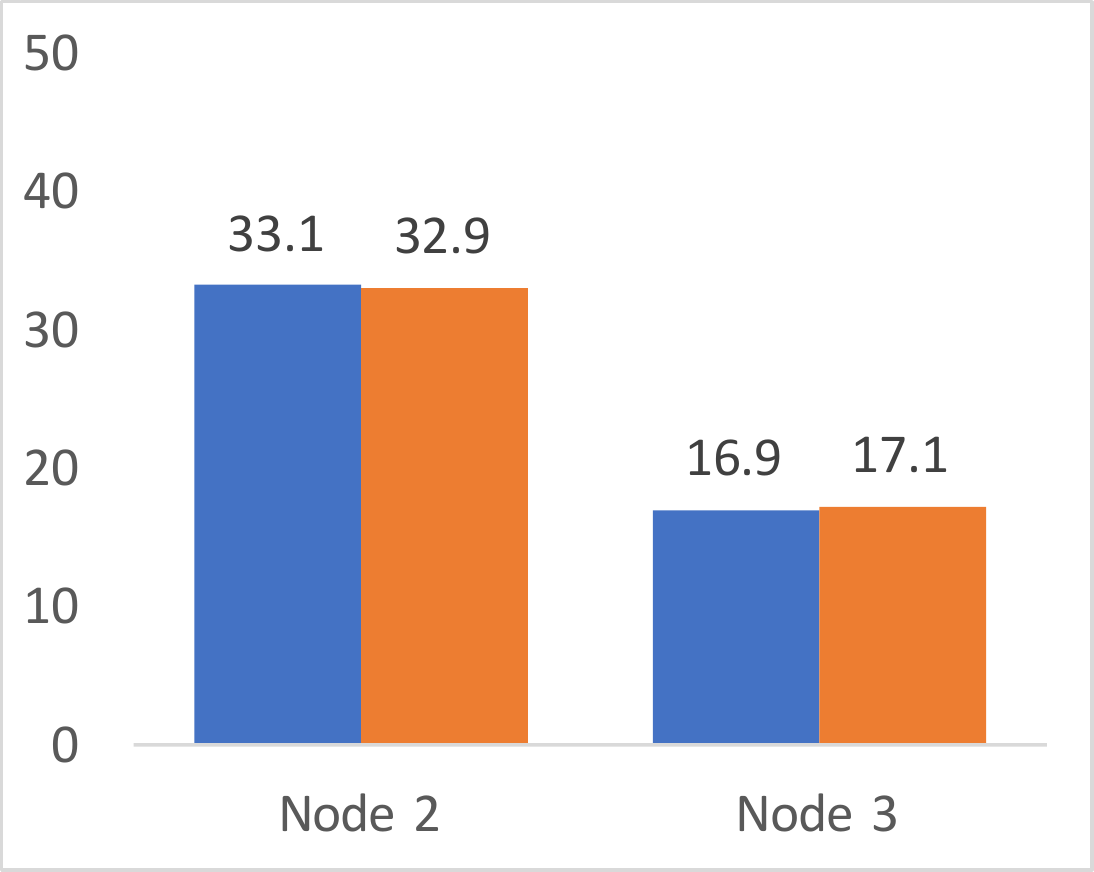}
			\label{fig:Driver_change}
		\end{subfigure}
		\begin{subfigure}[t]{.333\linewidth}
			\centering
			\caption{Riders Demand}
			\includegraphics[width=1\linewidth]{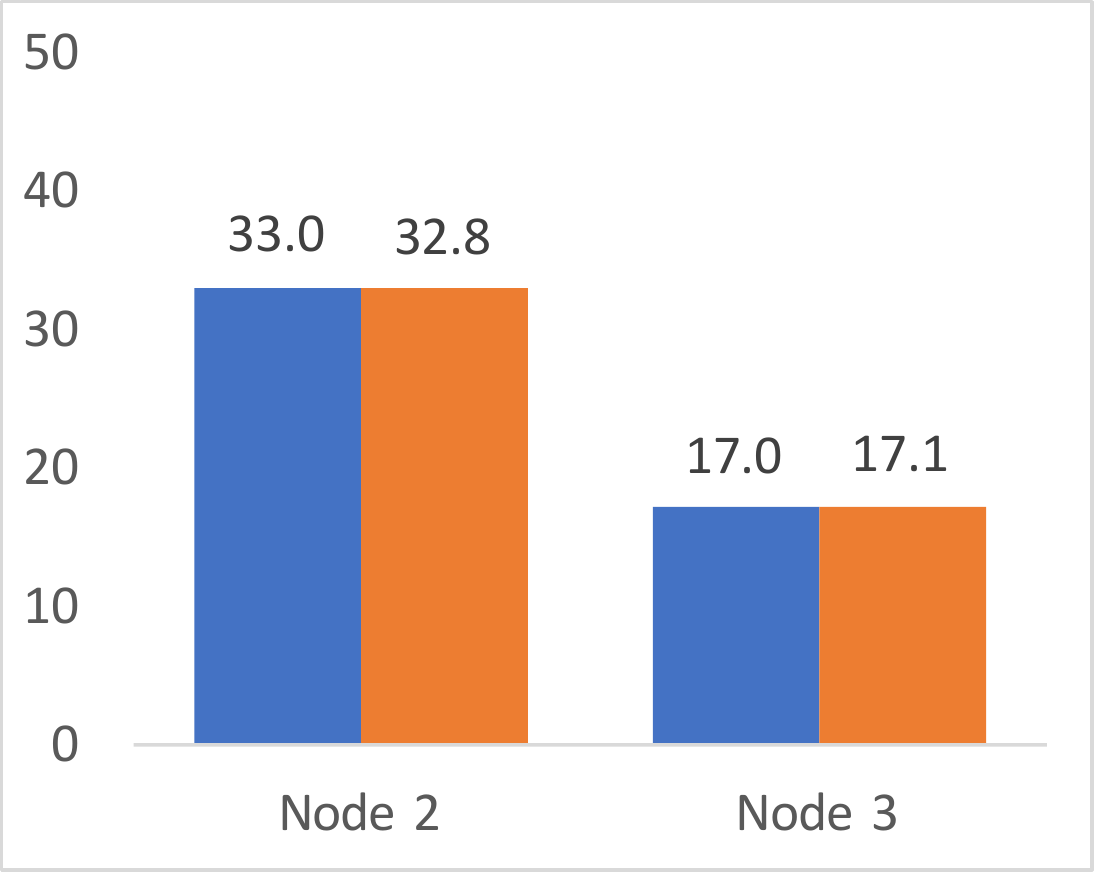}
			\label{fig:Rider_change}
		\end{subfigure}
		\begin{subfigure}[t]{.333\linewidth}
			\centering
			\caption{Equilbirium Prices $\bm{\rho}$}
			\includegraphics[width=1\linewidth]{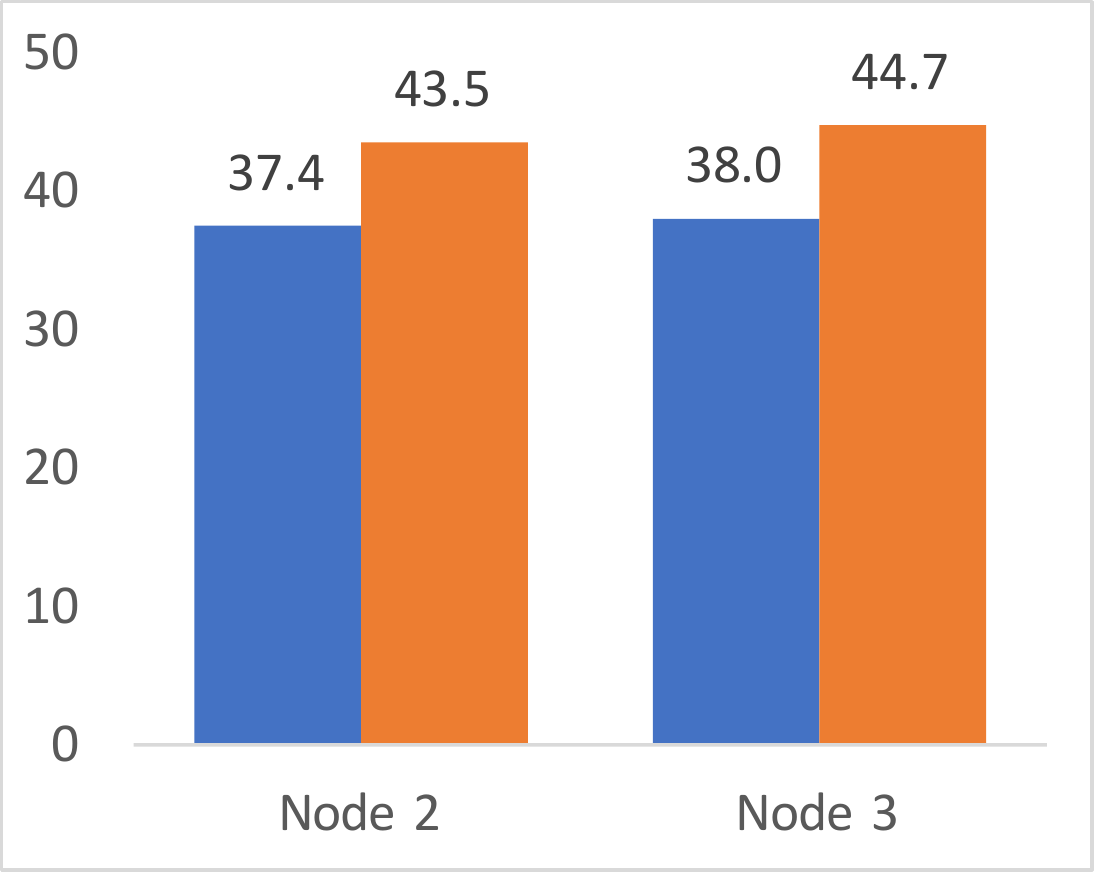} 
			\label{fig:Price_change}
	\end{subfigure}}
	\begin{subfigure}[t]{.533\linewidth}
			\centering
			\vspace{-0.5cm}
			\includegraphics[width=1\linewidth]{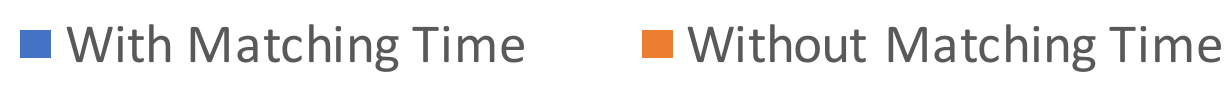} 
			\label{fig:legend_compare}
	\end{subfigure}
	\vspace{-2em}
	\caption{Equilibrium Solutions For Ignoring and Considering Matching Time Cases}
	\label{fig:ignore_consider_three_node}
\end{figure*}

We note that ignoring matching time is a special case of considering matching time when $a_0 = 0$ in (\ref{eq:waiting_time_flow}). We compare the equilibrium solutions between ignoring and considering matching time cases, with the results of equilibrium spatial pricing, driver supply and rider demand shown in Figure \ref{fig:ignore_consider_three_node}. We can see that, compared to equilibrium prices, driver and rider flows are less sensitive to whether matching time is considered or not. The reason is that different locations, when supply and demand are balanced, have similar waiting time. For drivers, their relocation choices only depend on the utility difference between locations, which will not change if similar waiting time is ignored. Therefore, driver supply remain unchanged (see Figure \ref{fig:Driver_change}). Given similar amount of driver supply, TNC will need to raise up locational prices (see Figure \ref{fig:Price_change}) to compensate the ignored waiting time, so that rider demand still balance driver supply (see Figure \ref{fig:Rider_change}).

\subsection{Sioux Falls Network}

\begin{figure}[htbp]
\begin{center}
    \includegraphics[width=0.5\textwidth]{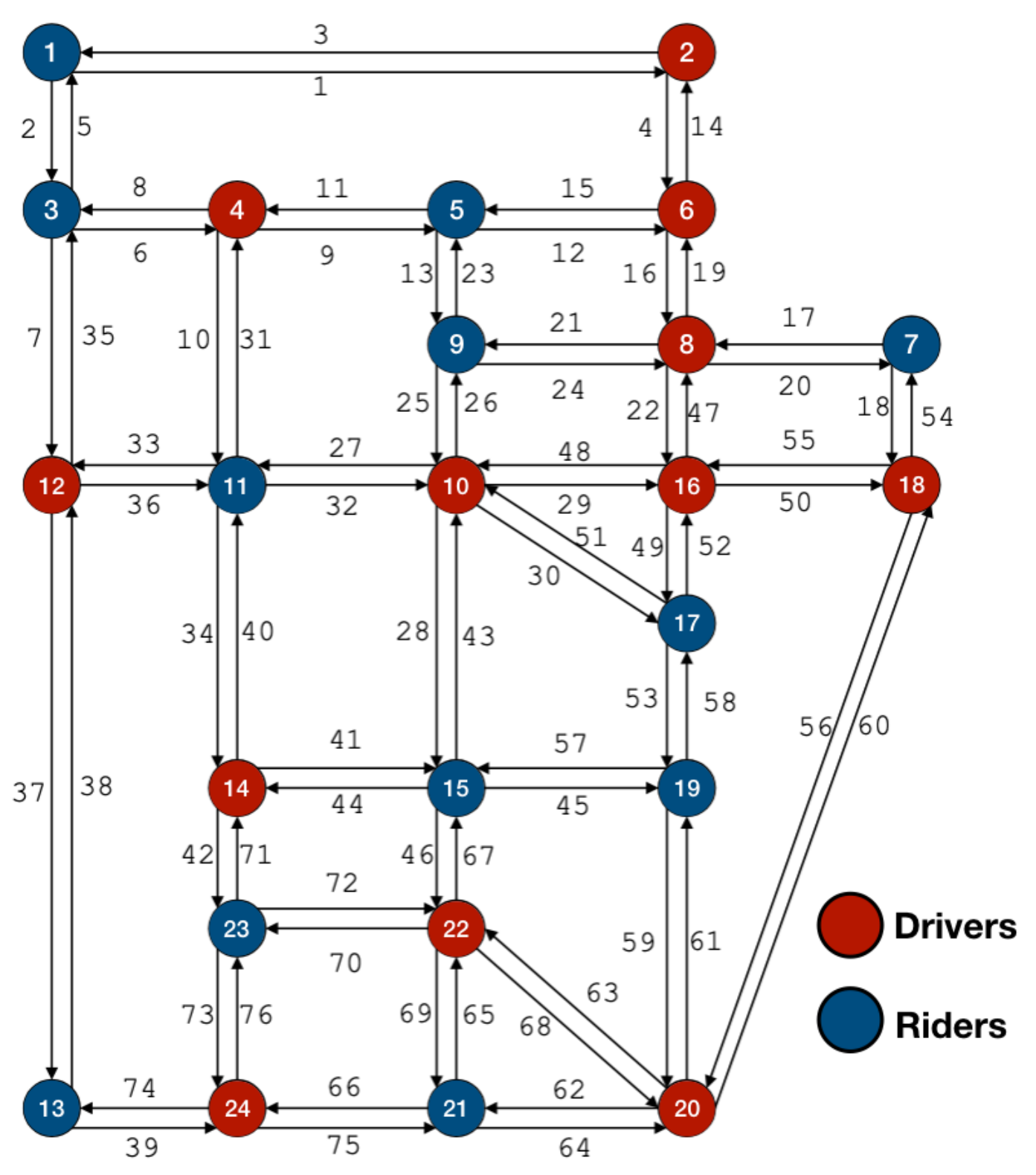}
\caption{Sioux Falls Test Network}
\label{fig:siou_fall}
\end{center}
\end{figure}

We implement our model using a medium-size test system, Sioux Falls network\footnote{Sioux Falls network consists of 24 nodes and 76 directed links. The number on each node/link in Figure \ref{fig:siou_fall} is the node/link index. Figure \ref{fig:siou_fall}  is edited base on a map provided by Hai Yang and Meng Qiang, Hong Kong University of Science and Technology.}, as shown in Figure \ref{fig:siou_fall}. In Figure \ref{fig:siou_fall}, the red and blue nodes (12 of each) represent the  original locations of drivers and riders, respectively. Drivers supply at each red node is 50, while the demand function at each blue node is $d_s = 300 - 5\rho_s$. We adopt a $4^{th}$-order Bureau of Public Roads (BPR) function for link cost: $t_a = t_a^0[1+0.15*(v_a/c_a)^4]$. Coefficients of utility function for base case (\ref{eq:utility}) are $\beta_0 = 0, \beta_1 = 1, \beta_2 =0.6$.

\subsubsection{Ignoring Matching Time}
\begin{figure*}[htbp]
	\centering
	\makebox[1.0\linewidth][c]{
		\begin{subfigure}[t]{.5\linewidth}
			\centering
			\caption{Locational Prices}
			\includegraphics[width=1\linewidth]{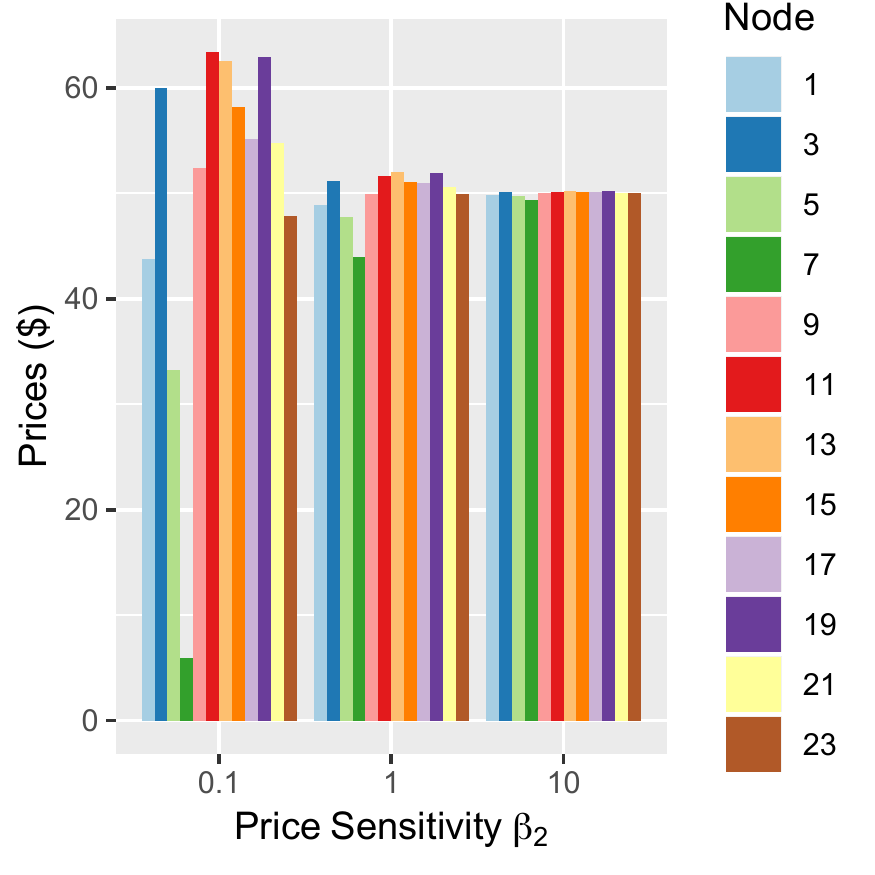}
			\label{fig:price_sioux_fall}
		\end{subfigure}
		\begin{subfigure}[t]{.5\linewidth}
			\centering
			\caption{Total Travel Time}
			\includegraphics[width=1\linewidth]{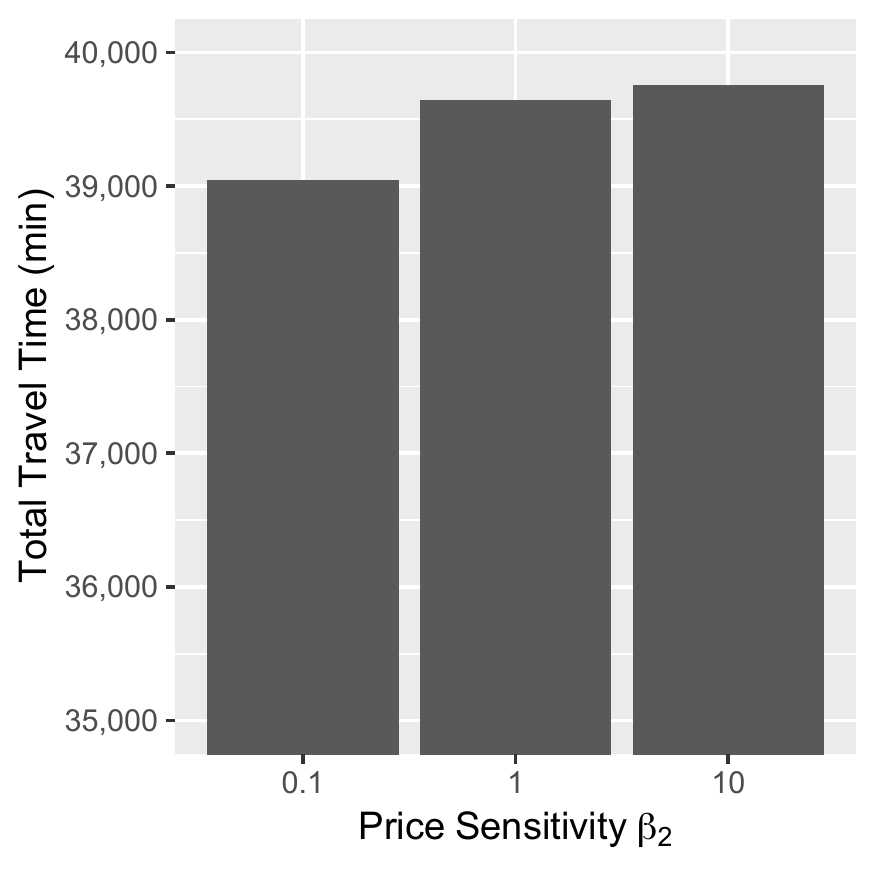} 
			\label{fig:time_sioux_fall}
	\end{subfigure}}
	\vspace{-2em}
	\caption{Impacts of Supply and Demand on Optimal Surge Pricing}
	\label{fig:senstivity_b2}
\end{figure*}

The optimal locational pricing decisions and total equilibrium travel time are shown in Figure \ref{fig:senstivity_b2}. With price sensitivity coefficient $\beta_2$ increase from 0.1 to 10, the price variances decrease (see Figure \ref{fig:price_sioux_fall}). This is because when drivers are more sensitive to prices, TNC needs smaller price difference to attract drivers to the desired locations. At the same time, drivers are more willing to travel longer distance in order to reach a location with higher price. This will lead to higher system total travel time, as shown in Figure \ref{fig:time_sioux_fall}. 

SCIP is not able to solve the MINLP reformulation of problem (\ref{model:upper}) for Sioux Falls network, while IPOPT can efficiently solve (\ref{model:combine}) in 11.3, 6.1, and 4.9 seconds for $b_2$ = 0.1, 1, and 10, respectively. 

\subsubsection{Considering Matching Time}

\begin{figure*}[htbp]
	\centering
	\makebox[1.0\linewidth][c]{
		\begin{subfigure}[t]{.5\linewidth}
			\centering
			\caption{Convergence of ES}
			\includegraphics[width=1\linewidth]{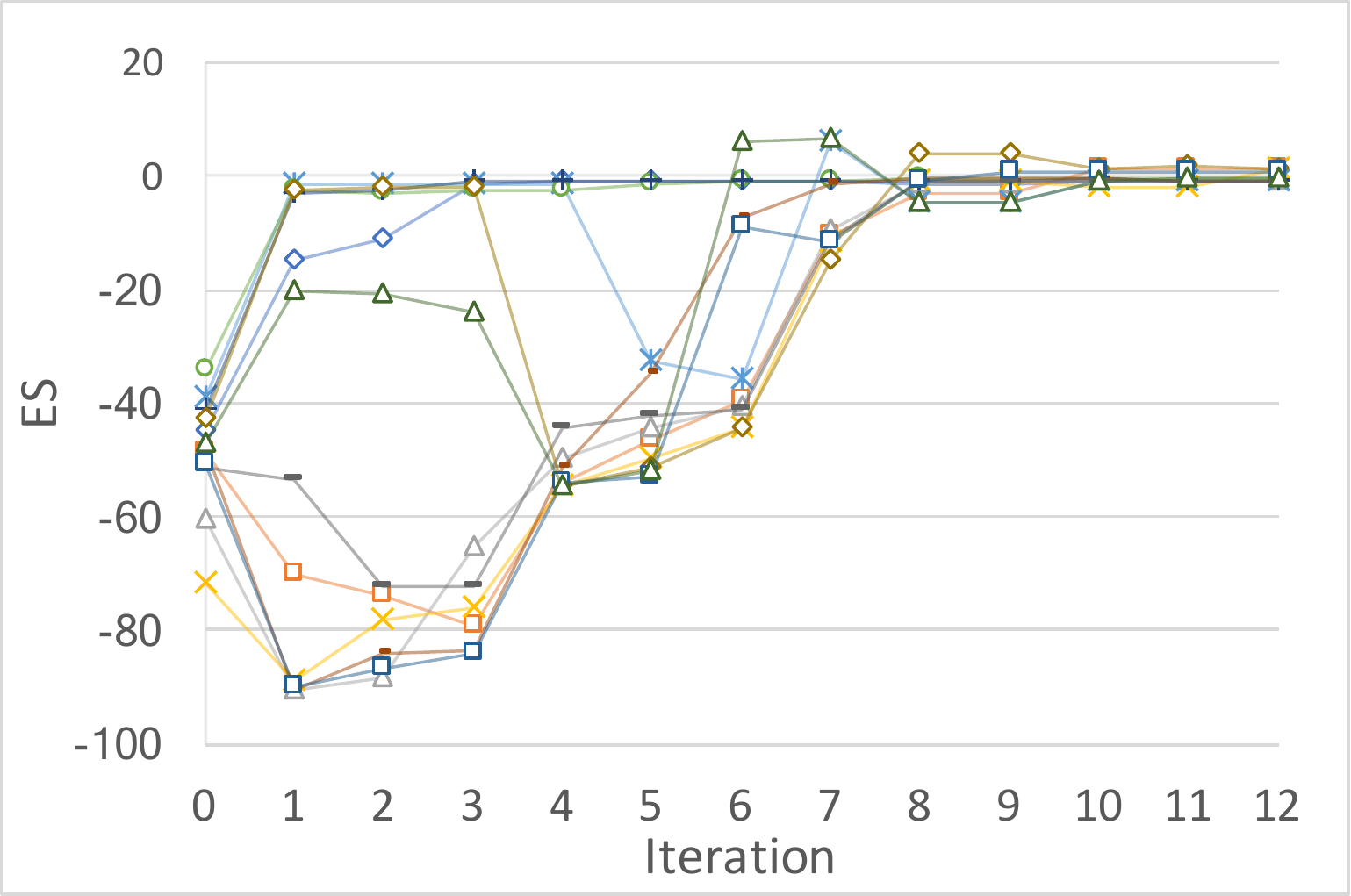}
			\label{fig:sioux_fall_matching_es}
		\end{subfigure}
		\begin{subfigure}[t]{.5\linewidth}
			\centering
			\caption{Convergence of $\bm{\rho}$}
			\includegraphics[width=1\linewidth]{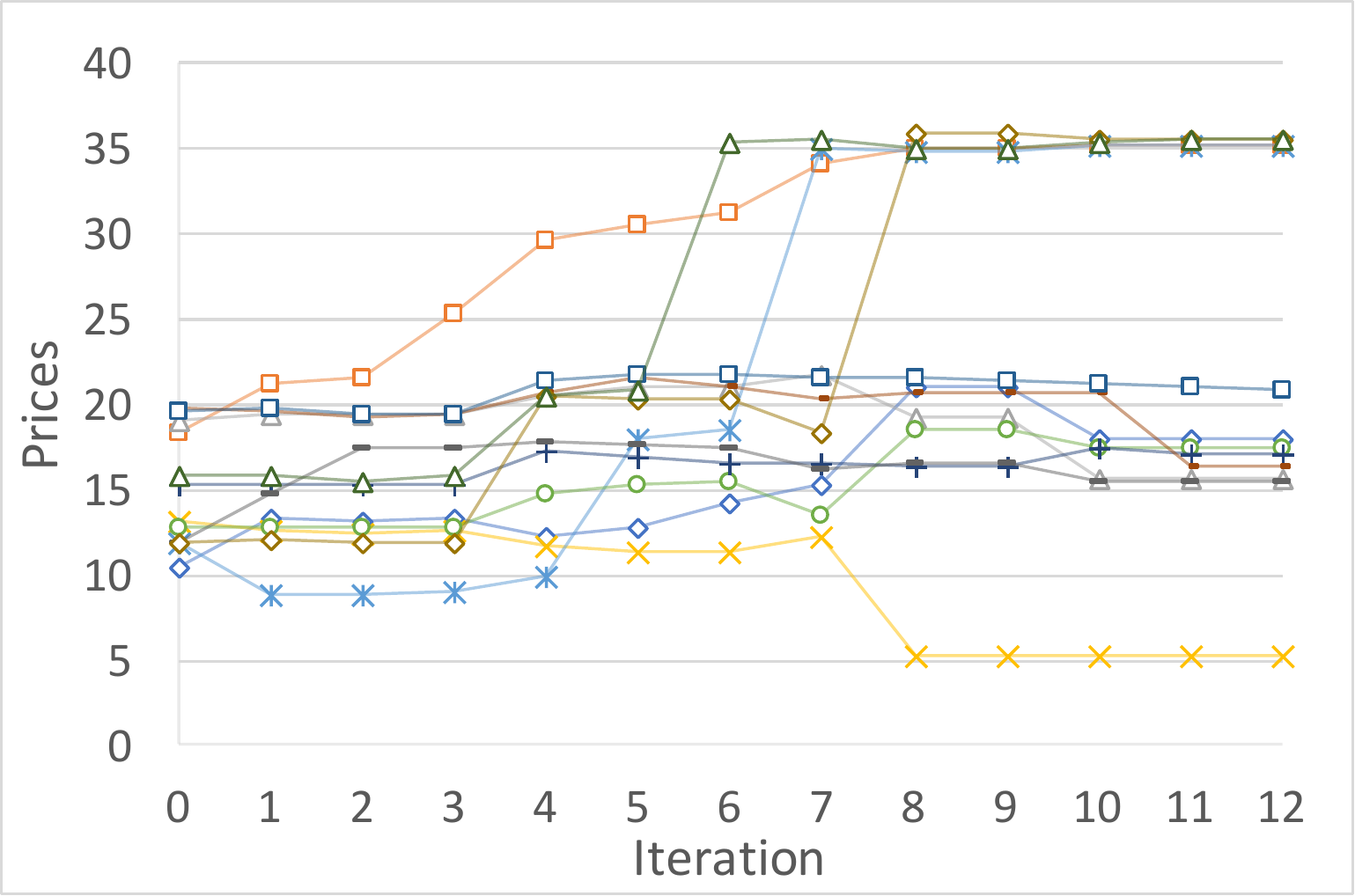} 
			\label{fig:sioux_fall_matching_prices}
	\end{subfigure}}
	\begin{subfigure}[t]{0.8\linewidth}
			\centering
			\vspace{-0.5cm}
			\includegraphics[width=1\linewidth]{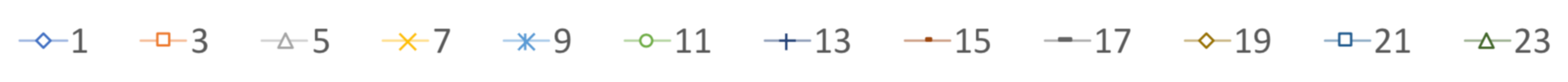} 
			\label{fig:legend_compare}
	\end{subfigure}
	\vspace{-2em}
	\caption{Convergence of ES and $\bm{\rho}$ Considering Matching Time}
	\label{fig:siou_fall_matching_conv}
\end{figure*}

Considering matching time, the problem can be solved by Algorithm \ref{alg:iter} in 6.4 hours, which is significantly longer than the computation time for ignoring matching time case. This is because of the lack of convexity for both upper and lower problems in (\ref{model:upper_matching}). But the algorithm converge reliably in 12 iterations, with convergence patterns of ES and $\bm{\rho}$ shown in Figure \ref{fig:siou_fall_matching_conv}.

Comparing between the impact of spatial pricing on transportation congestion, the total travel cost is 3805.5 when spatial pricing is adopted. When TNC adopts uniform pricing, the total travel cost is 3770.4. This observation is again consistent with the results in ignoring matching time case.

\section{Discussion}

In this paper, we have presented a new modeling framework and effective computational techniques for ride-sourcing spatial pricing problems, considering leader-follower structure of TNCs, riders, and drivers. In the lower level, interactions between drivers’ relocation, riders' mode choice, all travelers' routing behaviors, and network congestion are explicitly modeled. In the upper level, a monopoly TNC optimally determines spatial prices to achieve its own objectives. When matching time can be ignored, we show that the problem is equivalent to a convex optimization problem, which can be efficiently solved to global optimal by commercial nonlinear solver, such as IPOPT. Existence and uniqueness of optimal solutions are proved. When matching time cannot be ignored, we show that the lower-level drivers' and riders' two-sided interaction can be reformulated in a unified CDA framework with non-separable link-cost function through network augmentation. The problem can be solved reliably by approximating maxinf-points algorithm. 

In addition to showing the effectiveness of our solution approaches, we demonstrated that (1) transportation congestion affects spatial pricing strategies, and spatial pricing strategies may aggravate transportation congestion compared to uniform pricing; (2) the spatial pricing strategies minimizing supply-demand imbalance also maximize TNC profits when maximum demand and supply are comparable; (3) drivers' preference will significantly impact spatial prices and transportation congestion; and (4) the impacts of waiting time may be small for the distribution of drivers and riders, but may be large for the scale of spatial pricing.

This research can be extended in several directions. From a methodological viewpoint, our model mainly focuses on spatial aspects. Coupling temporal and spatial dimensions are critical to fully understand the impacts of dynamic pricing on transportation systems. In addition, several assumptions in this study can be relaxed. For example, the assumption on monopoly competition and inelastic total driver supply could be relaxed to generalize our findings. After gaining a clear understanding about the impacts of saptial pricing on transportation systems, strategies to mitigate the negative impacts from regulators perspective is urgently needed. Extending our model to a stochastic environment is another immediate next step considering the high level of uncertainties involved in ride-sourcing market, including congestion, rider demand, driver supply, etc. Although our model can be efficiently solved for small and medium transportation networks, algorithm development is needed for extreme large cases. Classic traffic assignment algorithm could be leveraged to solve the lower-level problem more efficiently. Solving a stochastic version of our model might also present further numerical challenges, for which decomposition methods (such as scenario decomposition and Benders decomposition) developed for stochastic programs may be integrated with the reformulation techniques in this study. From a pratical perspective, given the challenges of obtaining ride-sourcing data, agent-based simulation techniques can be leveraged to simulate the ride-sourcing environment and generate simulated data for implementation, validation and numerical insights.

\section*{References}
\bibliography{ref.bib}
\pagebreak 
\appendix
%
\input{./Appendics/app_pfs}

\pagebreak 
\setcounter{table}{0}
\input{./Appendics/app_explain}

\end{document}

%% file: mac.tex
\usepackage[usenames,dvipsnames]{color}
\usepackage{lineno,hyperref}
\usepackage{amsmath}
\usepackage{breqn}
\usepackage{calc}  
\usepackage{enumitem}  
\usepackage{amsfonts}
\usepackage{mathrsfs}
\usepackage{algpseudocode}
\usepackage{multirow}
\usepackage{bigstrut}
\usepackage{booktabs}
\usepackage{rotating}
\usepackage{hyperref}
\usepackage[top=3.5cm,bottom=3.5cm,left=3.5cm,right=3.5cm]{geometry}
\usepackage{caption}
\usepackage{subcaption}

\newtheorem{lemma}{Lemma}
\newtheorem{theorem}{Theorem}

\newtheorem{defn}{Definition}

\newcommand{\reals}{{I\kern-.35em R}}
\newcommand{\nats}{{I\kern-.35em N}}
\newcommand{\upto}{{\raise 1pt \hbox{$\scriptstyle \,\nearrow\,$}}}
\newcommand{\downto}{{\raise 1pt \hbox{$\scriptstyle \,\searrow\,$}}}
\newcommand{\eop}
	{\hfill{$\vcenter{\hrule height1pt \hbox{\vrule width1pt height5pt 
   	 \kern5pt \vrule width1pt} \hrule height1pt}$} \medskip}
\def\state #1. { \noindent{\bf#1.\enspace}}
\newcommand{\nargmaxinf}{\mathop{\rm argmaxinf}\nolimits}
\newcommand{\nargmin}{\mathop{\rm argmin}\nolimits}
\newcommand{\nargmax}{\mathop{\rm argmax}\nolimits}
\newcommand{\nsupinf}{\mathop{\rm supinf}\nolimits}
\newcommand{\ninf}{\mathop{\rm inf}\nolimits}

%% file: Appendics/app_pfs.tex
\section{Proofs.}\label{app:pfs}

\state Proof (Lemma \ref{lem:CDA}).
The proof of Lemma \ref{lem:CDA} is outlined as follows. Since (\ref{obj:lower}) is strictly convex (assuming $t_a(\cdot)$ is a strictly monotone increasing function) and constraints of model (\ref{model:lower}) are linear, Karush-Kuhn-Tucker (KKT) conditions are the sufficient and necessary optimality conditions. One can show that the KKT conditions are equivalent to the Wardrop equilibrium conditions (i.e., travel time in all routes actually used are equal and less than those that would be experienced by a single vehicle on any unused route) and logit model solutions. For details, please refer to \citep{evans1976derivation}.
\eop

\state Proof (Lemma \ref{lem:balance_matches}).
Denote the optimized drivers distribution of model (\ref{model:lower}) as $q_{rs}^{\ast}$. Because of constraint (\ref{cons:q_Q}), the total drivers supply is:
\[\sum_{r\in R}\sum_{s \in S}q_{rs}^{\ast} = \sum_{r\in R}Q_{r} \doteq \bar{Q},\] which is a constant.

Because of constraint (\ref{cons:matches}), $ m_s \leq \sum_{r \in \mathcal{R}} q_{rs}^{\ast}$. Therefore, $\sum_{s \in \mathcal{S}} m_s \leq \sum_{s \in \mathcal{S}}\sum_{r \in \mathcal{R}} q_{rs}^{\ast} = \bar{Q}$.  Denote the optimal value of $m_s$ as $m^{\ast}_s$, we have: 
\begin{equation}
    \sum_{s \in \mathcal{S}}m^{\ast}_s \leq \bar{Q}\label{eq:matches_cap}
\end{equation}

On the other hand, 
\begin{equation}
    \sum_{s \in \mathcal{S}}m^{\ast}_s \geq 0\label{eq:matches_cap_lower}
\end{equation}
The equality in (\ref{eq:matches_cap_lower}) holds if and only if 
\begin{equation}
    \sum_{r\in \mathcal{R}}q_{rs}^{\ast} = D_s - b_s\rho_s, \forall s \in \mathcal{S}\label{eq:balance}
\end{equation}
In other words, we just need to show that $\exists \bm{\rho} \in \reals^{S}$, s.t. (\ref{eq:balance}) holds.

Define excess supply $ES^s (\bm{\rho}) = \sum_{r \in \mathcal{R}}q_{rs}^{\ast} (\bm{\rho}) - (D_s - b_s \rho_s), \forall s$. For each location $s$, the ride-sourcing supply $\sum_{r \in \mathcal{R}}q_{rs}^{\ast} (\bm{\rho}) \leq \bar{Q}$; the equality holds when all the drivers relocate to location $s$. Therefore, the minimum locational price $\underline{\rho}^s = (D^s-\bar{Q})/b^s$. On the other hand, the minimum ride-sourcing supply is strictly positive. Therefore, the maximum locational price $\bar{\rho}_s = D_s/b_s$. This implies that $\rho_s \in [\underline{\rho}_s, \bar{\rho}_s], \forall s$.

Define a convex compact set $H = \Pi_{s\in \mathcal{S}}[\underline{\rho}_s, \overline{\rho}_s]$. Next, we construt a function $\boldsymbol{\rho}' = h(\boldsymbol{\rho}) = (h_s(\rho_s), \forall s)$, such that:
\[ \rho'_s = h_s(\rho_s) = \rho_s - \frac{(\overline{\rho}_s - \underline{\rho}_s) ES_s(\boldsymbol{\rho})}{\bar{Q}} = \rho_s - \frac{ES_s(\boldsymbol{\rho})}{b_s}, \rho_s \in [\underline{\rho}_s, \overline{\rho}_s] \]

We would like to show that $\rho'_s \in [\underline{\rho}_s, \overline{\rho}_s]$. Firstly, it can be easily seen that $h_k(\underline{\rho}_s) = \underline{\rho}_s - \frac{ES_s(\boldsymbol{\rho})}{b_s} = \underline{\rho}_s - \frac{\sum_{r \in \mathcal{R}}q_{rs}^{\ast} (\bm{\rho}) - (D_s - b_s \underline{\rho}_s)}{b_s} \leq \frac{D_s-\bar{Q}}{b_s}- \frac{0 - \bar{Q}}{b_s}= \overline{\rho}_s $, and $h_k(\overline{\rho}_s) = \overline{\rho}_s - \frac{ES_s(\boldsymbol{\rho})}{b_s} = \overline{\rho}_s - \frac{\sum_{r \in \mathcal{R}}q_{rs}^{\ast} (\bm{\rho}) - (D_s - b_s \overline{\rho}_s)}{b_s} \geq \frac{D_s}{b_s} - \frac{(Q - 0)}{b_s} = \underline{\rho}_s$. Secondly, we only need to show $h_s(\rho_s)$ is a monotone decreasing function. To show that, we take first order derivative of $h_s(\rho_s)$ with respect to $\rho_s$:
\[ \dot{h}_s(\rho_s) = 1 - \frac{1}{b_s}\frac{\partial ES_s(\boldsymbol{\rho})}{\partial\rho_s } = 1 - \frac{1}{b_s}(\frac{\partial \sum_{r \in \mathcal{R}}q_{rs}^{\ast} (\bm{\rho})}{\partial\rho_s } - \frac{\partial(D_s - b_s \rho_s)}{\partial\rho_s })\]

Because $\frac{\partial \sum_{r \in \mathcal{R}}q^{rs\ast} (\bm{\rho})}{\partial\rho^s } \geq 0$, 
\[ \dot{h}_s(\rho^s) \leq 1 - \frac{1}{b^s} (0 + b^s) = 0\]
 
Therefore, we have proved function $h(\cdot)$ is a mapping from $H$ to $H$. 
In addition, because $h(\cdot)$ is monotone and continuous; and $H$ is a convex compact subset of Euclidean space, based on Brouwer's fixed-point theorem, there exists $\boldsymbol{\rho}^{\ast} \in H$, such that $\boldsymbol{\rho}^{\ast} = h(\boldsymbol{\rho}^{\ast})$. Therefore, for any $s$, $\rho_s^{\ast} = \rho_s^{\ast} - \frac{ES_s(\boldsymbol{\rho}^{\ast})}{b^s}$, which leads to $ES_s(\boldsymbol{\rho}^{\ast}) = 0, \forall s$.

To show the uniqueness of $\bm{\rho}^{\ast}$, assuming that there exist $\bm{\rho}^{\ast\ast}$ such that $ES_s(\boldsymbol{\rho}^{\ast\ast}) = 0, \forall s$. Therefore,  $h(\boldsymbol{\rho}^{\ast\ast}) = \boldsymbol{\rho}^{\ast\ast}$, i.e., $\boldsymbol{\rho}^{\ast\ast}$ is another fixed point of $\boldsymbol{\rho}' = h(\boldsymbol{\rho})$. Without lost of generality, assuming ${\rho}_s^{\ast\ast} > {\rho}_s^{\ast}$. Since $h_s(\cdot)$ is monotone decreasing, $\rho_s^{\ast\ast} = h_s(\rho_s^{\ast\ast}) \leq h_s(\rho_s^{\ast}) = \rho_s^{\ast}$, which is a contraction with the assumption that ${\rho}_s^{\ast\ast} > {\rho}_s^{\ast}$.
\eop

\state Proof (Theorem \ref{thm:combine}).
Because problem (\ref{model:combine}), (\ref{model:lower}) and (\ref{model:demand}) are convex optimization, KKT conditions are both necessary and sufficient conditions. It can be shown that the KKT conditions of problem (\ref{model:combine}) are identical to the KKT conditions of problem (\ref{model:lower}), (\ref{model:demand}) and equilibirum constraint (\ref{model:equilibrium_constraints}). Therefore, $\bm{\rho}^{\ast}$ that solves problem (\ref{model:combine}) also soves MOPEC problem (\ref{model:lower}), (\ref{model:demand}) and (\ref{model:equilibrium_constraints}). Based on Lemma \ref{lem:balance_matches},  $\bm{\rho}^{\ast}$ solves original problem (\ref{model:upper}).

\eop

\state Proof (Lemma \ref{lem:mieq}).
On one hand, if $\rho^*$ is a maxinf-point of the Walrasian with $W(\rho^*,\cdot) \geq 0$, 
it follows that for all unit vectors $e^s=(0,\ldots,1,\ldots)$, the $s$-th 
entry is 1, $W(\rho^*,e^s) \geq 0$ which implies ${\rm ES}_s(\rho^*) = 0$. On the other hand, if $\rho^*$ is an equilibrium prices, i.e. ${\rm ES}_s(\rho^*) = 0$, it follows that $W(\rho^*,\cdot) = 0$. Because $W(\rho,\cdot) \leq 0$ for any $\rho$, $\rho^*$ is a maxinf point of the Walrasian.

\eop

%% file: Appendics/app_explain.tex
\section{Matching Time $\sim$ Driver and Rider Flow}\label{app:explain}

\begin{figure}[htbp]
\begin{center}
    \includegraphics[width=0.5\textwidth]{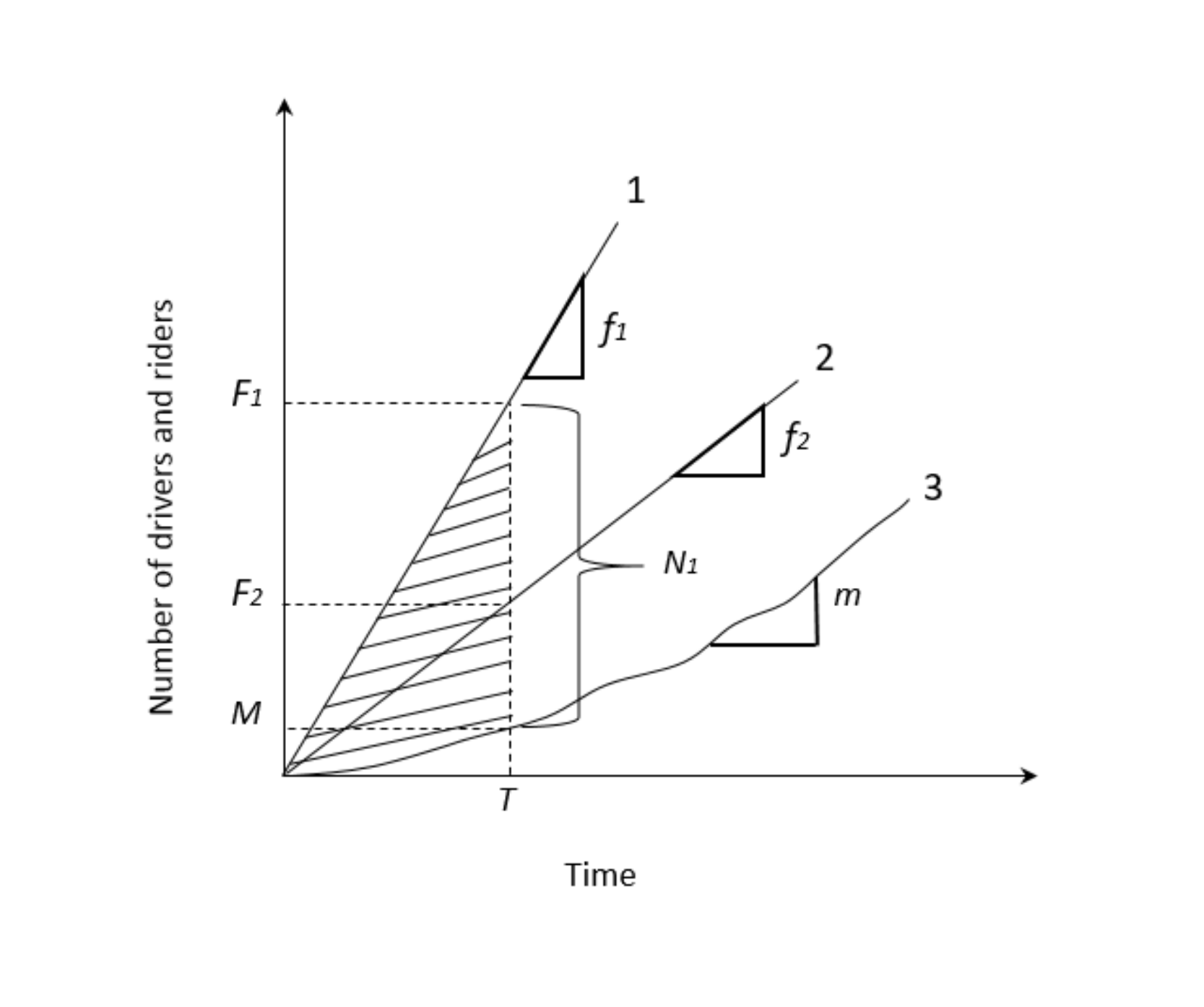}
    \vspace{-1cm}
\caption{Graphical Representation of Matching Time}
\label{fig:match_time}
\end{center}
\end{figure}

Equation (\ref{eq:waiting}) can be better explained graphically in Figure \ref{fig:match_time}. Line 1 and 2 indicates the accumulative drivers (slope $f_{1}$) and riders (slope $f_{2}$) arrival curves, respectively. Line 3 represents the accumulative matching (with slope $m$ calculated from Cobb-Douglas matching function). The area of the shadow areas indicate the total delay of drivers (or riders) till time period $T$, from which, the average delay per driver ($t_1$) and per rider ($t_2$) can be calculated. By varying $f_1$ and $f_2$, we are able to calculate $t_1$ and $t_2$ numerically, which will be used for estimating parameters in (\ref{eq:waiting_time_flow}). 

\begin{figure}[htbp]
\begin{center}
    \includegraphics[width=0.7\textwidth]{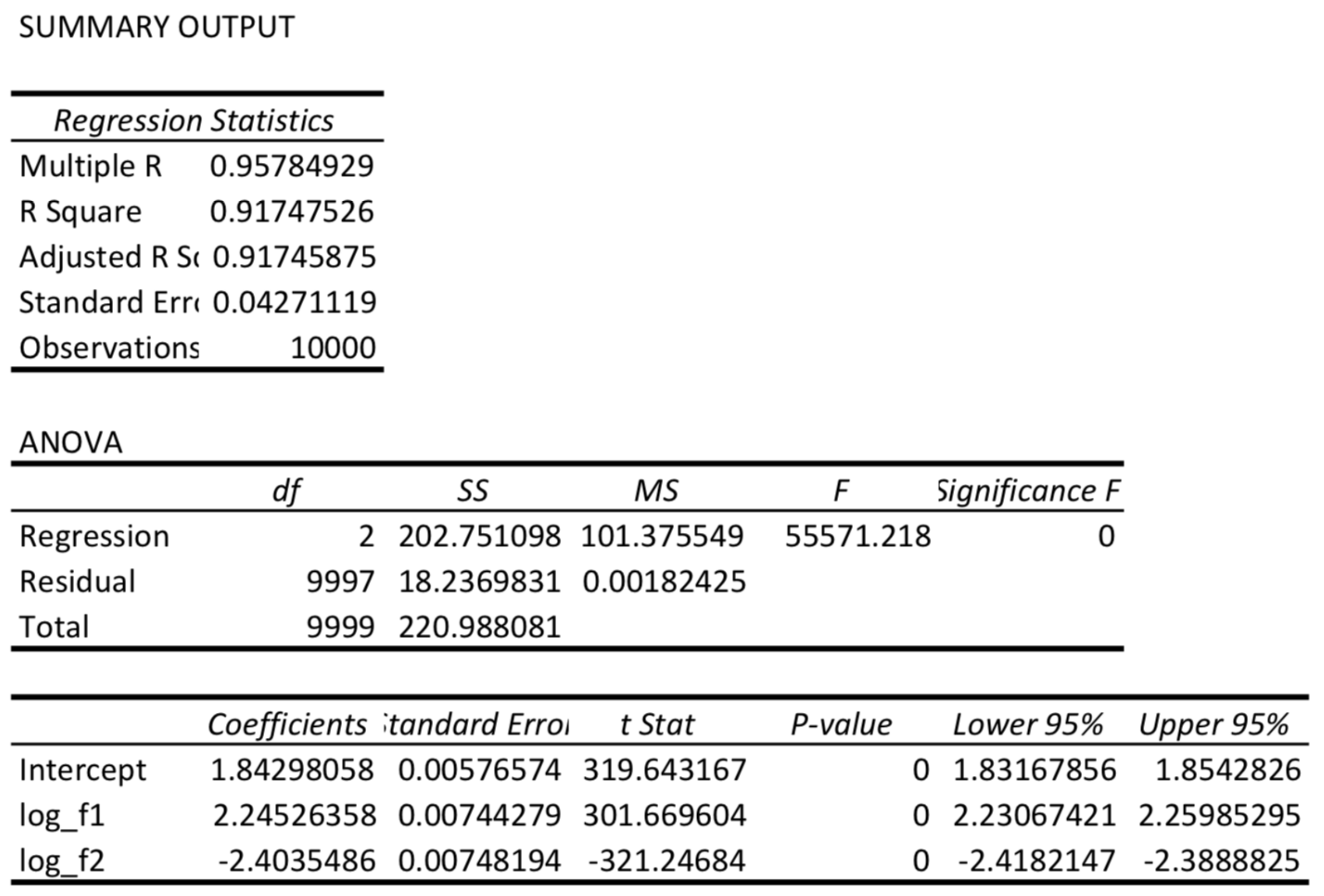}
\caption{Regression Summary Output of (\ref{eq:waiting_time_flow})}
\label{fig:regression_output}
\end{center}
\end{figure}